\documentclass[a4paper,10pt]{article}

\usepackage[title]{appendix}
\usepackage{setspace}
\usepackage[bottom]{footmisc}
\usepackage{color}
\usepackage{subfig}
\usepackage{svg}
\usepackage{hyperref}
\usepackage{multirow}
\usepackage{cite}
\usepackage{makecell}
\usepackage{amsmath} 
\usepackage{mathtools}
\usepackage{cleveref}
\usepackage{graphicx}
\usepackage{epsfig}
\usepackage{epstopdf}
\usepackage{float}
\usepackage{titling}
\usepackage[margin=1.0in]{geometry}
\usepackage{etoolbox}

\title{An analytical solution for supersonic flow over a circular cylinder using an optimized shock shape}
\makeatletter
\renewcommand\@date{{%
  \vspace{-\baselineskip}%
  \large\centering
  \begin{tabular}{@{}c@{}}
    S R Siva Prasad Kochi\textsuperscript{1} \\
    \normalsize siva.ksr@gmail.com
  \end{tabular}%
  \quad and\quad
  \begin{tabular}{@{}c@{}}
    M Ramakrishna\textsuperscript{2} \\
    \normalsize krishna@ae.iitm.ac.in
  \end{tabular}

  \bigskip

  \textsuperscript{1}Dept. of Aerospace Engg., IIT Madras.\par
  \textsuperscript{2}Professor, Dept. of Aerospace Engg., IIT Madras.

  \bigskip

  \today
}}
\makeatother
\begin{document}

\maketitle

\begin{abstract}
 An analytical solution for high supersonic flow over a circular cylinder based on Schneider's inverse method \cite{schneider} has been presented. In the inverse method, a shock shape is assumed and the corresponding flow field and the shape of the body producing the shock are found by integrating the equations of motion using the stream function. A shock shape theorised by Moeckel \cite{moeckel} has been assumed and it is optimized by minimising the error between the shape of the body obtained using Schneider's method and the actual shape of the body. A further improvement in the shock shape is also found by using the Moeckel's shock shape in a small series expansion. With this shock shape, the whole flow field in the shock layer has been calculated using Schneider's method by integrating the equations of motion. This solution is compared against a fifth order accurate numerical solution using the discontinuous Galerkin method (DGM) and the maximum error in density is found to be of the order of $10^{-3}$ which demonstrates the accuracy of the method used for both plane and axisymmetric flows.
 
 {{\bf Keywords:} cylinder, supersonic flow, shock wave, shock shape, inverse method, discontinuous Galerkin method}
\end{abstract}

\section{Introduction}\label{intro}

Analytical methods for hypersonic flows past blunt bodies in plane and axisymmetric flows have been investigated quite extensively in the past. These analyses can be found in papers like those of Chester \cite{chester1}, \cite{chester2} and Freeman \cite{freeman1} and in the books by Hayes and Probstein \cite{hprobstein} and Rasmussen \cite{rasmussen1}. There are two different problems that are examined by them. The first one is the so called direct problem where the shape of the body is given and the shock shape and the flow field are to be found. The second one is the inverse problem where the shape of the shock wave is given, and it is required to find the corresponding flow field and the shape of the body producing the shock. We are interested in the inverse problem and particularly the method devised by Schneider in \cite{schneider} for flow past blunt bodies to solve the problem. We use the shock shape theorised by Moeckel \cite{moeckel} for blunt bodies along with the method of Schneider \cite{schneider} to find an accurate analytical solution for plane flow past a circular cylinder for high supersonic velocities. In a way, this procedure is similar to the shock fitting technique \cite{salas1} used by numerical solvers where the shock shape is assumed and improved further based on the solution. 
\\
\\
\noindent The analysis developed by Schneider allows a very elegant treatment of the inviscid hypersonic blunt body problem. The fundamental advantage of this method is its applicability in the flow field from the stagnation region up to large distances from the nose of the body. This method for the solution in shock layer is based on two main assumptions. First, it is assumed that the density immediately behind a strong shock is much larger than in front of the shock. The second assumption is based on the pressure along a streamline. With these two assumptions, Schneider obtained the solution by integrating the equations of motion using the stream function. The flow properties immediately behind the shock are obtained using the Rankine-Hugoniot equations. As a result of these assumptions and the subsequent solution, it is not necessary to have a thin shock layer for this solution to hold. Due to these advantages, this method has been used extensively in various aerospace applications like in \cite{chou1}, \cite{maslen1}, and \cite{nagaraja1}. Recently, this method has also been used for astrophysical calculations in \cite{sb1} and \cite{sfkwbw}. This method was also extended to three-dimensional flows by Schwarze in \cite{schwarze1}. Though this method was developed for hypersonic flows, we show that this works well even in the high supersonic flow regime.
\\
\\
\noindent Our method of solution is as follows. We assume a shock shape in the form given by Moeckel in \cite{moeckel} as a single parameter family of hyperbolas. We determine that parameter by minimising the error between the shape of the body obtained using Schneider's method and the actual shape of the body. This gives a better shock shape $f$. Now, we improve this shock shape further by writing it as a series expansion ($af+bf^{2}$) with unknown coefficients. We determine these coefficients using the same optimization and get a much better approximation for the shock shape. We use this shock shape to find the solution of the flow in the shock layer and compare the solution obtained with a numerical solution calculated using a fifth order accurate discontinuous Galerkin method (DGM) and overset grids where the shock has been captured within a grid line \cite{srspkmr4}. We find that the maximum error in density is of the order of $10^{-3}$ (for the shock shape $af+bf^{2}$) demonstrating the accuracy of the solution method. We also found that using more terms in the series expansion for the shock shape (like $af+bf^{2}+cf^{3}$, $af+bf^{2}+cf^{3}+df^{4}$) does not improve the accuracy of the solution.
\\
\\
\noindent The paper is organized as follows. We describe the formulation of the Schneider's method used for all our results in Section \ref{formulation}, approximation of the shock shape is described in Section \ref{approxShockShape}, the results obtained are given in Section \ref{results} and we conclude the paper in Section \ref{conc}. We also give the formulation and validation of the numerical method used (discontinuous Galerkin method) in Appendix \ref{app:DGMFormulationValidation}.

\section{Formulation}\label{formulation}

\begin{figure}[htbp]
\begin{center}
\scalebox{1.5}{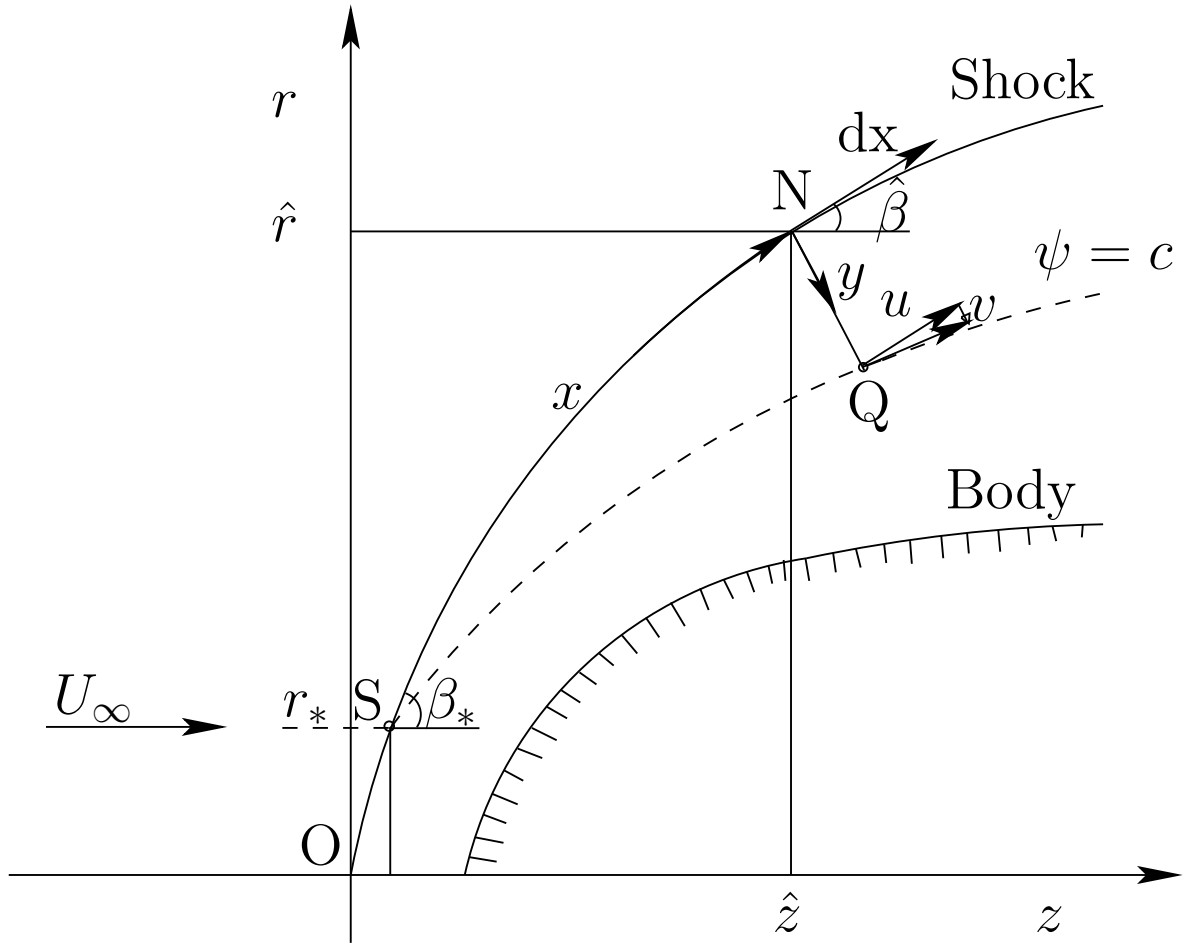}
\caption{The shock, body and a streamline are shown with a shock oriented coordinate system. Here, $z$ and $r$ are the Cartesian coordinates for plane flow with origin O, $x$ and $y$ are the distances along the shock surface and normal to it and Q is an arbitrary point in the shock layer. $u$ and $v$ are the velocities in $x$ and $y$ directions. $U_{\infty}$ is the freestream velocity. $\beta$ is the shock angle. $\psi$ is the stream function. Flow quantities immediately behind the shock at the point N (point where a normal to the shock intersects Q) are denoted by a hat ($\wedge$), and at the point S (intersection of the streamline through Q with the shock) by an asterisk (*). This diagram is adapted from Schneider's paper \cite{schneider}.}
\label{fig:ShockBody}
\end{center}
\end{figure}

\noindent We now describe the formulation of Schneider's method in brief. We have mostly followed the notation used in \cite{sb1} for the description of this method for the plane or axisymmetric supersonic flow past a blunt body. The shape of the bow shock wave (and hence the shock angle $\beta$) is given, and it is required to find the corresponding flow field and the shape of the body producing the shock. Let the flow be inviscid and without heat conduction, and let the gas be perfect. $z$ and $r$ are the Cartesian coordinates for plane flow, and also the cylindrical coordinates for axisymmetric flow. We also use a shock oriented curvilinear coordinate system, where $x$ is the distance along the shock surface in the plane formed by the shock normal and the direction of the uniform fluid flow, and $y$ is the distance normal to the shock surface as shown in Figure \ref{fig:ShockBody}. Consider an arbitrary point Q in the shock layer. A streamline $\psi=c$ passes through this point. The streamline intersects the shock at the point S. The point N on the shock is determined such that a normal to the shock at N also intersects the streamline, $\psi=c$, at the point Q. We can now determine the $x$ and $y$ coordinates of Q and the corresponding velocity components $u$ and $v$.
The $z$-axis is taken to be parallel to the direction of the incident flow. We emphasize that this is only assumed for convenience, but is not a general restriction. It can be deduced from Figure \ref{fig:ShockBody} that

\begin{equation}\label{coord:1}
 z = \hat{z} + y \sin \hat{\beta}
\end{equation}

\begin{equation}\label{coord:2}
 r = \hat{r} - y \cos \hat{\beta}
\end{equation}

\noindent where $\hat{\beta}$ is the shock inclination angle at the point N ($\hat{z}$, $r$). The flow quantities immediately behind the shock at the point N are denoted by a hat ($\wedge$), and at the point S by an asterisk (*). Undisturbed flow quantities far upstream are denoted by the subscript $\infty$. 
\\

\noindent Since the functions $\hat{z}(x)$ and $\hat{r}(x)$ are known for a given shock shape, equations \eqref{coord:1} and \eqref{coord:2} can be used to calculate the coordinates $z$ and $r$ of a point Q from its coordinates $x$ and $y$. The curvature of the shock contour at the point N is denoted by $\hat{\kappa}(x)$, defined as positive when the surface is concave on the side of positive $y$.
\\
\noindent Now, the continuity equation can be written as

\begin{equation}\label{mass}
 \frac{\partial (r^{j} \rho u)}{\partial x} + \frac{\partial [(1-\hat{\kappa}y) r^{j} \rho v]}{\partial y} = 0
\end{equation}

\noindent where $\rho$ is the fluid density. The parameter j is 0 for plane flow and 1 for axisymmetric flow. We now define a stream function $\psi$ which satisfies \eqref{mass} by

\begin{equation}\label{stream1}
 \frac{\partial \psi}{\partial x} = (1-\hat{\kappa}y)r^{j}\rho v,
\end{equation}

\begin{equation}\label{stream2}
 \frac{\partial \psi}{\partial y} = -r^{j}\rho u
\end{equation}

\noindent The stream function defined by \eqref{stream1} and \eqref{stream2} contains a constant of integration. This constant is chosen such that $\psi=0$ is the body stream line. Then the stream function represents the mass flow between the streamline $\psi$ = constant. and the surface of the projectile, per unit depth for plane flows, and per unit azimuthal angle (in radians) for axisymmetric flows. At the point N, the stream function therefore is
\\
\begin{equation}\label{streamP}
 \hat{\psi}_{N} = \rho_{\infty}U_{\infty}\frac{\hat{r}^{1+j}}{1+j}
\end{equation}

\noindent Also, $\psi$ is connected to the coordinate $r_{*}$ of the point S by 

\begin{equation}\label{streamS}
 \psi_{*} = \rho_{\infty}U_{\infty}\frac{r_{*}^{1+j}}{1+j}
\end{equation}

\noindent We now introduce a new coordinate system with $\psi$ and $\bar{x}=x$ as variables. With these new coordinates, the equations of energy, entropy, and momentum conservation can written as:

\begin{equation}\label{newEnergy}
 u^{2}+v^{2}+2h = u_{*}^{2} + v_{*}^{2} + 2h_{*} = \text{constant}.,
\end{equation}

\begin{equation}\label{newEntropy}
 \frac{\partial S}{\partial \bar{x}} = 0 \quad \text{or} \quad S = S_{*}(\psi)
\end{equation}

\begin{equation}\label{newXMomentum}
 u\frac{\partial u}{\partial \bar{x}} + v\frac{\partial v}{\partial \bar{x}} + \frac{1}{\rho}\frac{\partial P}{\partial \bar{x}} = 0 
\end{equation}

\begin{equation}\label{newYMomentum}
 (1-\hat{\kappa}y)(1-j\frac{y}{\hat{r}}\cos \hat{\beta})\hat{r}^{j}\frac{\partial P}{\partial \psi} = \hat{\kappa}u + \frac{\partial v}{\partial \bar{x}}
\end{equation}

\noindent where $P$ is the fluid pressure, $S$ the entropy, and $h$ is the specific enthalpy.
\\
\noindent The variable $y$ is now a dependent variable and we now get the transformation equations using \eqref{stream1} and \eqref{stream2} as

\begin{equation}\label{yDependent1}
 \frac{\partial y}{\partial \bar{x}} = (1-\hat{\kappa}y)\frac{v}{u}
\end{equation}

\begin{equation}\label{yDependent2}
 \frac{\partial y}{\partial \psi} = -\frac{1}{[1-j(y/\hat{r})\cos \hat{\beta}]\hat{r}^{j}\rho u}
\end{equation}

\noindent Now, the equations \eqref{newXMomentum}, \eqref{newYMomentum}, \eqref{yDependent1}, \eqref{yDependent2}, and \eqref{newEnergy} (or \eqref{newEntropy}) provide a set of five equations for the
five unknown dependent variables $u$, $v$, $y$, $P$ and $h$ (or $S$).
\\
\\
\noindent The flow quantities immediately behind the shock may be obtained from the Rankine-Hugoniot jump conditions in terms of the inverse compression ratio across the shock $\hat{\chi}=\rho_{\infty}/\hat{\rho}$. They are

\begin{equation}\label{rh1}
 \hat{u} = U_{\infty}\cos \hat{\beta}
\end{equation}

\begin{equation}\label{rh2}
 \hat{v} = U_{\infty} \hat{\chi} \sin \hat{\beta}
\end{equation}

\begin{equation}\label{rh3}
 \hat{P} = P_{\infty} + \rho_{\infty}U_{\infty}^{2}(1-\hat{\chi})\sin^{2} \hat{\beta}
\end{equation}

\begin{equation}\label{rh4}
 \hat{h} = h_{\infty} + \frac{1}{2}U_{\infty}^{2}(1-\hat{\chi}^{2})\sin^{2} \hat{\beta}
\end{equation}

\noindent These equations maintain their validity even if the hats are replaced by asterisks. 
\\
\noindent Schneider's method for the solution in shock layer is based on two main assumptions. First, it is assumed that the density immediately behind a strong shock is much larger than in front of the shock; i.e.,

\begin{equation}\label{assump1}
 \hat{\chi}=\frac{\rho_{\infty}}{\rho}\ll 1 \quad \text{and} \quad \chi_{*}=\frac{\rho_{\infty}}{\rho_{*}}=O(\hat{\chi})
\end{equation}

\noindent Second, the pressure $P$ at the point Q of the disturbed flow field is not much smaller than the pressure $\hat{P}$ at the point N (see Figure \ref{fig:ShockBody}); i.e.,

\begin{equation}\label{assump2}
 \frac{\hat{P}}{P}=O(1) \quad (\text{on }x=\text{const., }y>0)
\end{equation}

\noindent The symbol $f(x)=O(g(x))$ means that $|f(x)|$ is not very large in comparison with $|g(x)|$. Using these assumptions, Schneider in \cite{schneider} obtained the pressure at an arbitrary point Q approximately as

\begin{equation}\label{SchneiderPressure}
 P = \hat{P}-\frac{\hat{\kappa}}{\hat{r}^{j}}\int_{\psi}^{\hat{\psi}} \left[u_{*}^{2}+2\left(h_{*}-h(\hat{P},S_{*})\right)\right]^{1/2}d\psi'
\end{equation}

\noindent Here, we can note that all quantities on the right hand side of equation \eqref{SchneiderPressure} are given by the boundary conditions at the shock or by the equation of state. The terms excluded in making the approximation are coming from the stagnation region, as well as the region near the stagnation region where $u\gg \hat{u}$ – are of the order of $\hat{\chi}$ and hence they contribute only negligibly to the integral in equation \eqref{SchneiderPressure}. Therefore, the whole state of the gas is known in the streamline coordinate system $(\bar{x}, \psi)$ after evaluating $S=S_{*}(\psi)$ and $P$ from equation \eqref{SchneiderPressure}. The actual location of the body in space can be determined by solving the differential equation \eqref{yDependent2}. We do this by separation of variables giving the distance from the shock surface $y$ as a function of $\bar{x}$ and $\psi$ as

\begin{equation}\label{yAndPsi}
 y\left(1-\frac{j\cos \hat{\beta}}{2\hat{r}}y\right) = \frac{1}{\hat{r}^{j}}\int_{\psi}^{\hat{\psi}} \frac{d\psi'}{\rho u}
\end{equation}

\noindent Neglecting errors of $O(\chi)$ and following \cite{hprobstein}, we may replace the velocity component $u$ by

\begin{equation}\label{xVelocity}
 u^{2} = u_{*}^{2} + 2[h_{*}-h(P,S_{*})]
\end{equation}

\noindent This is obtained from the energy equation \eqref{newEnergy}. The integral in equation \eqref{yAndPsi} can now be written as

\begin{equation}\label{YandPsi}
 Y = \int_{\psi}^{\hat{\psi}}\frac{d\psi'}{\rho(P,S_{*})[u_{*}^{2}+2(h_{*}-h(P,S_{*}))]^{1/2}}
\end{equation}

\noindent Solving the quadratic equation in $y$ on the left hand side of equation \eqref{yAndPsi}, we have to distinguish between plane ($j = 0$) and axisymmetric ($j = 1$) flows. This gives

\begin{equation}\label{planeY}
 \text{for } j=0: \quad y = Y;
\end{equation}

\begin{equation}\label{axisymmetricY}
 \text{for } j=1: \quad y = \frac{\hat{r}}{\cos \hat{\beta}}\left[1-\left(1-\frac{2Y\cos \hat{\beta}}{\hat{r}^{2}}\right)^{1/2}\right]
\end{equation}

\noindent Equations \eqref{YandPsi}, \eqref{planeY}, and \eqref{axisymmetricY} give us the required results in $\bar{x}$ and $\psi$ coordinates. We can now use the transformations \eqref{coord:1} and \eqref{coord:2} to get the required results in $z$ and $r$ coordinates.
\\
For a perfect gas with constant specific heats, the inverse compression ratio is given by

\begin{equation}\label{compressRatio}
 \chi_{*} = \frac{\gamma -1}{\gamma +1} + \frac{2}{(\gamma +1)M_{\infty}^{2}\sin^{2}\beta_{*}}
\end{equation}

\noindent We point out that an analogous relation is valid for $\hat{\chi}$, if all asterisks are replaced by hats in equation \eqref{compressRatio}.
\\
Using the shock conditions \eqref{rh1}-\eqref{rh4}, together with \eqref{compressRatio}, the two integrals \eqref{SchneiderPressure} and \eqref{YandPsi}, which have to be evaluated, become

\begin{equation}\label{finalP}
 P = \hat{P} - \frac{U_{\infty}\hat{\kappa}}{\hat{r}^{j}} \int_{\psi}^{\hat{\psi}} \left[\cos^{2} \beta_{*} + \left(\frac{2}{(\gamma -1)M_{\infty}^{2}}+\sin^{2}\beta_{*}\right)\left(1-\left(\frac{\sin^{2}\hat{\beta}}{\sin^{2}\beta_{*}}\right)^{\frac{\gamma -1}{\gamma}}\right)\right]^{1/2}d\psi'
\end{equation}

\begin{equation}\label{finalY}
 Y = \frac{1}{\rho_{\infty}U_{\infty}} \int_{\psi}^{\hat{\psi}} \chi_{*}\left(\frac{\hat{P}\sin^{2}\beta_{*}}{P\sin^{2}\hat{\beta}}\right)^{1/\gamma}\left[\cos^{2}\beta_{*}+\left(\frac{2}{(\gamma -1)M_{\infty}^{2}}+\sin^{2}\beta_{*}\right)\left(1-\left(\frac{P\sin^{2}\hat{\beta}}{\hat{P}\sin^{2}\beta_{*}}\right)^{\frac{\gamma -1}{\gamma}}\right)\right]^{-1/2}d\psi'
\end{equation}

\noindent The curvature of a curve in space can be calculated using the formula

\begin{equation}\label{curvatureFormula}
 \kappa = \frac{|d^2r/dz^{2}|}{\left[1+(dr/dz)^{2}\right]^{3/2}}
\end{equation}

\noindent The density is obtained from the formula

\begin{equation}\label{densityFormula}
 \rho = \rho_{*}\left(\frac{P}{P_{*}}\right)^{1/\gamma} = \left(\frac{\gamma -1}{\gamma +1}+\frac{2}{(\gamma +1)M_{\infty}^{2}\sin^{2}\beta_{*}^{2}}\right)^{-1}\rho_{\infty}\left(\frac{P\sin^{2}\hat{\beta}}{\hat{P}\sin^{2}\beta_{*}}\right)^{1/\gamma}
\end{equation}

\noindent For simplicity we introduce dimensionless units following \cite{sb1}. The normalized stream function then becomes $\Psi=\psi/\rho_{\infty}U_{\infty}L^{j+1}$, where $L$ is a characteristic length.
\\
\\
On the surface of the body we need $\psi=0$, so the pressure on the body surface $P_{b}(x)$, as well as the shock layer thickness can be obtained by replacing the lower limits in equations \eqref{finalP} and \eqref{finalY} by zero as:

\begin{equation}\label{bodyPressure}
 \frac{P_{b}}{\rho_{\infty}U_{\infty}^{2}} = \frac{1}{\gamma M_{\infty}^{2}} + (1-\hat{\chi})\sin^{2}\hat{\beta} - \frac{\hat{\kappa}}{\hat{r}^{j}} \int_{0}^{\hat{\Psi}} \left[\cos^{2}\beta_{*}+\left(\frac{2}{(\gamma -1)M_{\infty}^{2}}+\sin^{2}\beta_{*}\right)\left(1-\left(\frac{\sin^{2}\hat{\beta}}{\sin^{2}\beta_{*}}\right)^{\frac{\gamma -1}{\gamma}}\right)\right]^{1/2}d\Psi
\end{equation}

\noindent where the relation 

\begin{equation}\label{pressureRelation}
 \frac{\hat{P}}{\rho_{\infty}U_{\infty}^{2}} = \frac{P_{\infty}}{\rho_{\infty}U_{\infty}^{2}} + (1-\hat{\chi})\sin^{2}\hat{\beta} = \frac{1}{\gamma M_{\infty}^{2}} + (1-\hat{\chi})\sin^{2}\hat{\beta}
\end{equation}

\noindent has been used. Also, we get

\begin{equation}\label{shockOffsetPlane}
 \text{for } j=0: \quad \Delta = Y\Bigr|_{\psi=0}
\end{equation}

\begin{equation}\label{shockOffsetAxisymmetric}
 \text{for } j=1: \quad \Delta = \frac{\hat{r}}{\cos \hat{\beta}}\left[1-\left(1-\frac{2Y\Bigr|_{\psi=0}\cos \hat{\beta}}{\hat{r}^{2}}\right)^{1/2}\right]
\end{equation}

\noindent with

\begin{equation}\label{YatZero}
 Y\Bigr|_{\psi=0} = \int_{0}^{\hat{\Psi}}\chi_{*}\left(\frac{\hat{P}\sin^{2}\beta_{*}}{P\sin^{2}\hat{\beta}}\right)^{1/\gamma}\left[\cos^{2}\beta_{*}+\left(\frac{2}{(\gamma -1)M_{\infty}^{2}}+\sin^{2}\beta_{*}\right)\left(1-\left(\frac{P\sin^{2}\hat{\beta}}{\hat{P}\sin^{2}\beta_{*}}\right)^{\frac{\gamma -1}{\gamma}}\right)\right]^{-1/2}d\Psi
\end{equation}

\noindent As soon as the bow shock wave is parameterized, the required quantities can be computed by evaluating two integrals (namely equations \eqref{bodyPressure} and \eqref{shockOffsetPlane} or \eqref{shockOffsetAxisymmetric} and \eqref{YatZero}, together with the boundary values \eqref{newEntropy} and \eqref{rh1}-\eqref{rh4}). In deriving this result, Schneider \cite{schneider} emphasizes that it is \textit{not} necessary to have a thin shock layer for these equations to hold. Also, throughout the rest of the paper, we consider only plane flow with $j=0$.

\section{Approximation of shock shape}\label{approxShockShape}

\noindent We now need an appropriate shock shape to use the results of Section \ref{formulation}. For this, we use the approximate relation given by Moeckel in \cite{moeckel}. This has the desirable property of yielding a single relation between freestream Mach number and distance of the shock wave from the body. Detached shock waves are normal to the freestream direction on the axis of symmetry, and are asymptotic to the freestream Mach lines at great distances from the axis of symmetry. Hence it seems plausible that the shape of the wave may be approximated by a hyperbola. Moeckel therefore postulated that the equation of the shock is 
 
 \begin{equation}\label{eqnOfShock}
 \beta r = \sqrt{z^{2}-z_{0}^{2}}
 \end{equation}
 
 \noindent where $\beta = \sqrt{M_{\infty}^{2}-1}$ and $z_{0}$ is the location of the vertex of the wave. Transforming the coordinate system so that the origin is located on the shock, we get the shock shape to be

 \begin{equation}\label{eqnOfShockFinal}
 \beta r = \sqrt{(z+z_{0})^{2}-z_{0}^{2}}
 \end{equation}
 
\begin{figure}[htbp]
\begin{center}
\includegraphics[scale=1.0]{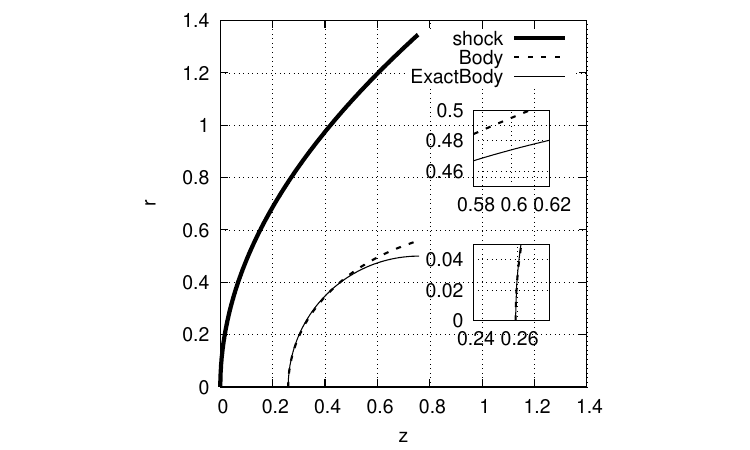}
\caption{The shock using equation \eqref{eqnOfShockFinal} with $z_{0}=17.615$ (thick solid line), the body obtained from \eqref{YatZero} (for plane flow) (thin dashed line), and the actual shape of the body (thin solid line) after minimising the error between the obtained shape of the body and the actual shape of the body for Mach number $M=4.0$. A zoomed-in view of the shock shape at two different locations ([0.24,0.275]$\times$[0,0.05], [0.58,0.62]$\times$[0.45,0.5]) is also shown in order to visualize the error clearly. Radius of the circular cylinder is $0.5$}
\label{fig:M4p0ShockBodyShapeOld}
\end{center}
\end{figure}

 \noindent Here $z_{0}$ is the only unknown and it represents the shock offset distance (which is the distance between the stagnation point on the body and the vertex of the shock wave). Starting from a reasonable approximation for this $z_{0}$, we can calculate the shape of the body using equation \eqref{YatZero} (for plane flow). By using the error between the obtained shape of the body and the actual shape of the body (in this case, one quarter of a circle), we optimize the parameter $z_{0}$ so that the error is minimised using secant method. After doing this, the shock, the body obtained from \eqref{YatZero} (for plane flow), and the actual shape of the body for Mach number $M=4.0$ are shown in Figure \ref{fig:M4p0ShockBodyShapeOld} where the radius of the circular cylinder is taken to be $0.5$.
 \\
 \\

 \noindent To get a better approximation for the shock shape, we use the relation given by Moeckel as follows. We first define
 
 \begin{equation}\label{fDef}
  f = \frac{1}{\beta}\sqrt{(z+z_{0})^{2}-z_{0}^{2}}
 \end{equation}
 
 \noindent Now, we approximate the shape of the shock again by using two more parameters $a$ and $b$ as
 
 \begin{equation}\label{shapeBetterApprox}
  r = af + bf^{2}
 \end{equation}

\begin{figure}[htbp]
\begin{center}
\includegraphics[scale=1.0]{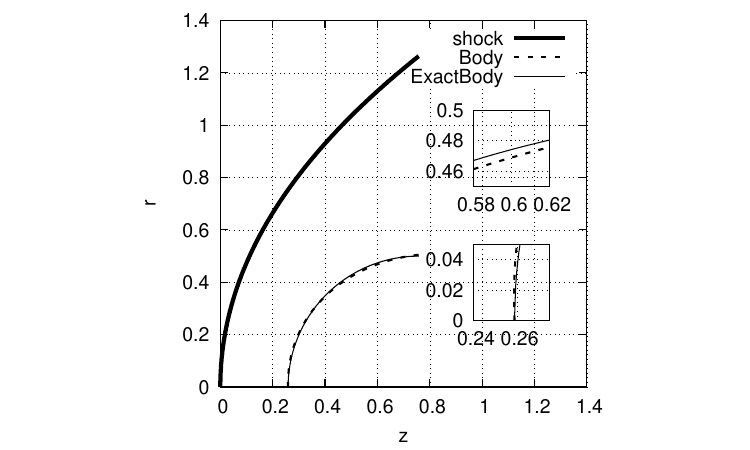}
\caption{The shock using equation \eqref{shapeBetterApprox} with $z_{0}=17.615$, $a=0.998$, and $b=-0.045$ (thick solid line), the body obtained from \eqref{YatZero} (for plane flow) (thin dashed line), and the actual shape of the body (thin solid line) after minimising the error between the obtained shape of the body and the actual shape of the body for Mach number $M=4.0$. A zoomed-in view of the shock shape at two different locations ([0.24,0.275]$\times$[0,0.05], [0.58,0.62]$\times$[0.45,0.5]) is also shown in order to visualize the error clearly. Radius of the circular cylinder is $0.5$}
\label{fig:M4p0ShockBodyShapeNew}
\end{center}
\end{figure}

 \noindent We now have three unknown parameters which are $z_{0}$, $a$ and $b$. We again determine these three parameters by optimizing them while minimising the error between the obtained shape of the body (using equation \eqref{YatZero}) and the actual shape of the body. After optimizing the three parameters $z_{0}$, $a$, and $b$, the shock, the body obtained from \eqref{YatZero} (for plane flow), and the actual shape of the body for Mach number $M=4.0$ are shown in Figure \ref{fig:M4p0ShockBodyShapeNew}. We can now clearly see the improvement over the previous result as the predicted body shape is very close to the actual shape of the body.
 \\

\begin{figure}[htbp]
\begin{center}
\includegraphics[scale=1.0]{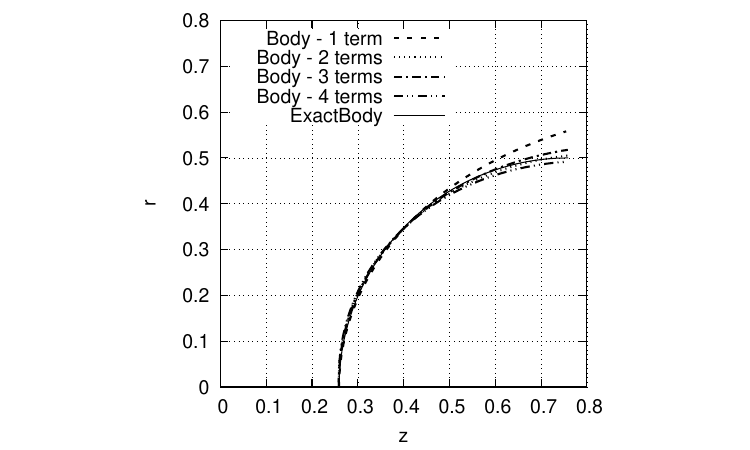}
\caption{The body obtained from \eqref{YatZero} (for plane flow) for shock shapes given by $f$ (dashed line), $af + bf^{2}$ (dotted line), $af + bf^{2} + cf^{3}$ (dash and dot line), $af + bf^{2} + cf^{3} + df^{4}$ (dash and double dot line) and the actual shape of the body (thin solid line) after minimising the error between the obtained shape of the body and the actual shape of the body for Mach number $M=4.0$. Radius of the circular cylinder is $0.5$}
\label{fig:M4p0ShockBodyShapeNewMoreTerms}
\end{center}
\end{figure}

\noindent We can continue this procedure with higher degree approximations in $f$ (i.e., $af + bf^{2} + cf^{3}$, $af + bf^{2} + cf^{3} + df^{4}$ and so on) by optimizing each of the coefficients of $f$ by minimising the error between the obtained shape of the body and the actual shape of the body. We show the results obtained for body shape using such approximations for shock shapes in Figure \ref{fig:M4p0ShockBodyShapeNewMoreTerms} for Mach number $M=4.0$. We can see from Figure \ref{fig:M4p0ShockBodyShapeNewMoreTerms} that using higher degree approximations (more than 2) does not really improve the body shape that much. Even the solution obtained using such higher degree approximations is not much better as will be explained in Section \ref{results}. For this reason, we have used equation \eqref{shapeBetterApprox} as the shock shape for all the remaining calculations unless otherwise specified.
\\

\noindent After obtaining the shock shape, we can now get the full flow field between the shock and the body on the streamlines by evaluating equations \eqref{finalP}, \eqref{finalY}, and \eqref{densityFormula}. We present all our results for flow over a circular cylinder in the next section.

\section{Results}\label{results}

\subsection{Plane flow:}\label{planeFlow}

We consider the plane flow past a circular cylinder for high supersonic Mach numbers and approximate the shape of the shock using equation \eqref{shapeBetterApprox} as discussed in Section \ref{approxShockShape}. We consider five different Mach numbers $M=4.0,5.0,6.0,7.0,$ and $8.0$ to illustrate our results. For all these Mach numbers, the parameters $z_{0}$, $a$, and $b$ are optimized by minimising the error between the obtained shape of the body (using equation \eqref{YatZero}) and the actual shape of the body. The parameters obtained for the considered Mach numbers are tabulated in Table \ref{table:1}.
\\
\begin{table}[htbp]
\centering
\begin{tabular}{|c|c|c|c|}
\hline
Mach Number & $z_{0}$ & a & b \\
\hline
4.0 & 17.615 & 0.998 & -0.045 \\
\hline
5.0 & 26.755 & 0.998 & -0.050 \\
\hline
6.0 & 38.495 & 0.998 & -0.052 \\
\hline
7.0 & 51.982 & 0.998 & -0.054 \\
\hline
8.0 & 67.984 & 0.998 & -0.058 \\
\hline
\end{tabular}
\caption{The optimized parameters $z_{0}$, $a$, and $b$ in \eqref{shapeBetterApprox} obtained by minimising the error between the obtained shape of the body (using equation \eqref{YatZero}) and the actual shape of the body.}
\label{table:1}
\end{table}

\begin{figure}[htbp]
\begin{center}
\includegraphics[scale=0.22]{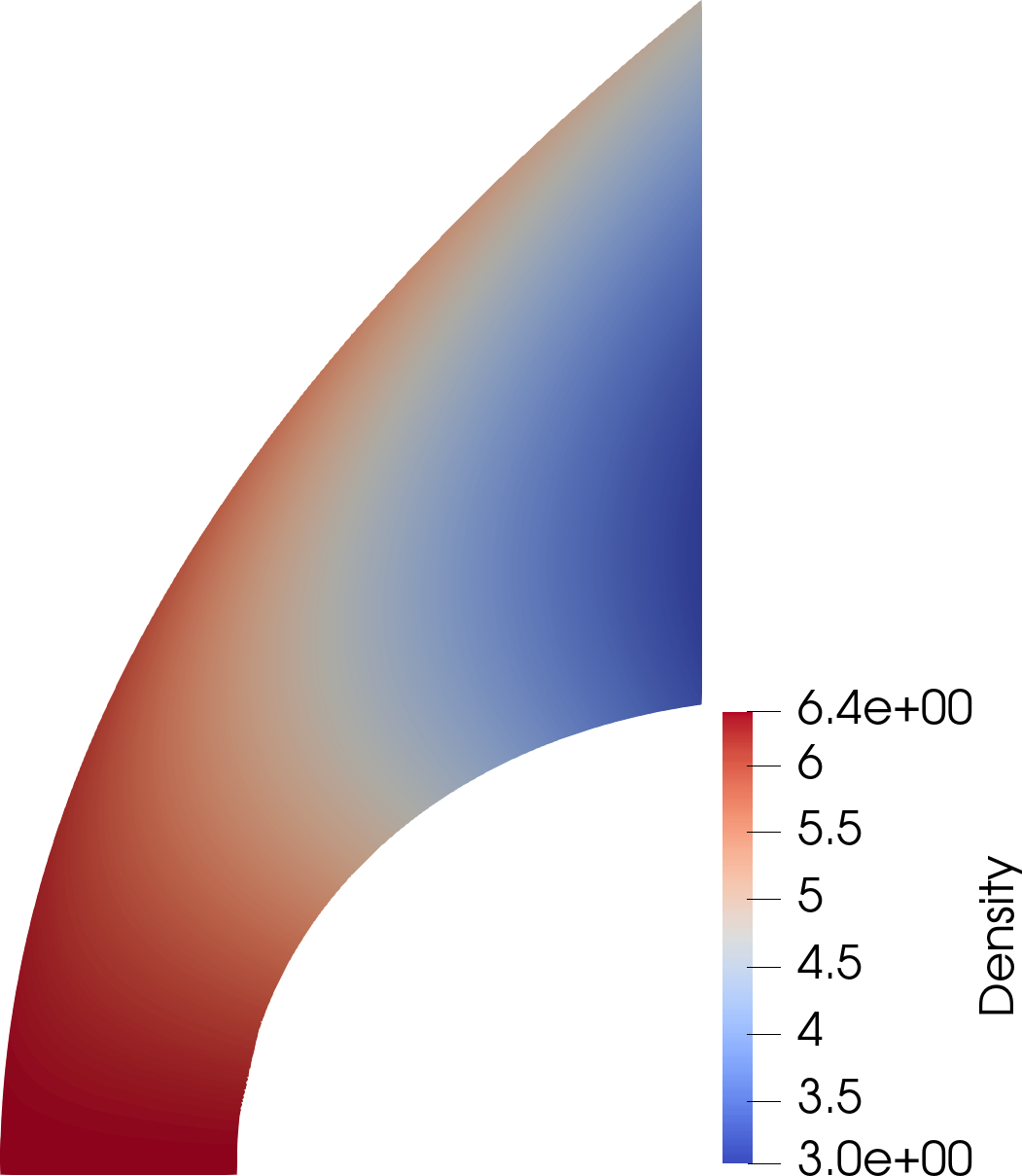}
\caption{Density variation for Mach number $4.0$ obtained using \eqref{finalP}, \eqref{finalY}, and \eqref{densityFormula} with shock shape given in \eqref{shapeBetterApprox} with parameters given in Table \ref{table:1} for 200 streamlines.}
\label{fig:M4p0Density}
\end{center}
\end{figure}

\begin{figure}[htbp]
\begin{center}
\includegraphics[scale=0.22]{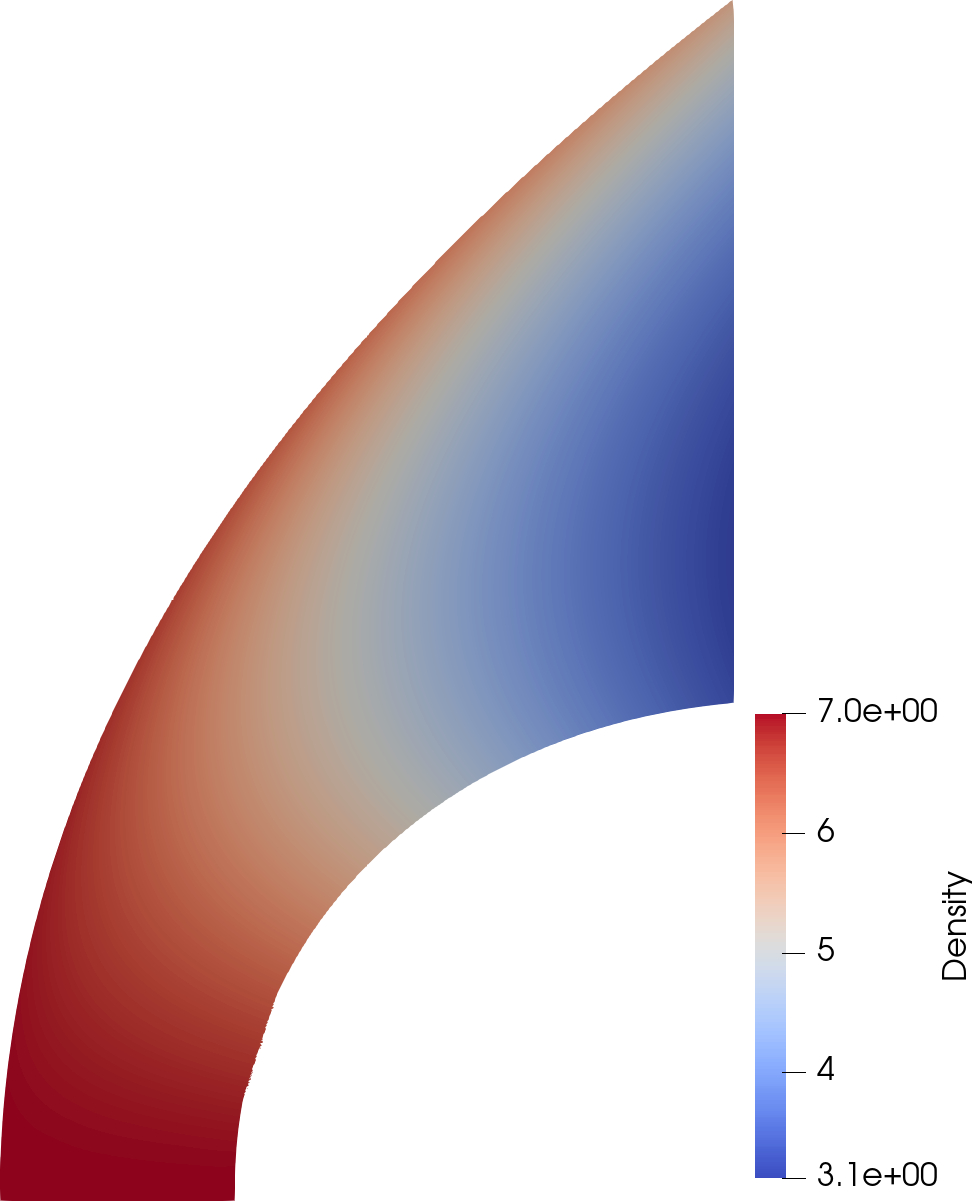}
\caption{Density variation for Mach number $5.0$ obtained using \eqref{finalP}, \eqref{finalY}, and \eqref{densityFormula} with shock shape given in \eqref{shapeBetterApprox} with parameters given in Table \ref{table:1} for 200 streamlines.}
\label{fig:M5p0Density}
\end{center}
\end{figure}

\begin{figure}[htbp]
\begin{center}
\includegraphics[scale=0.22]{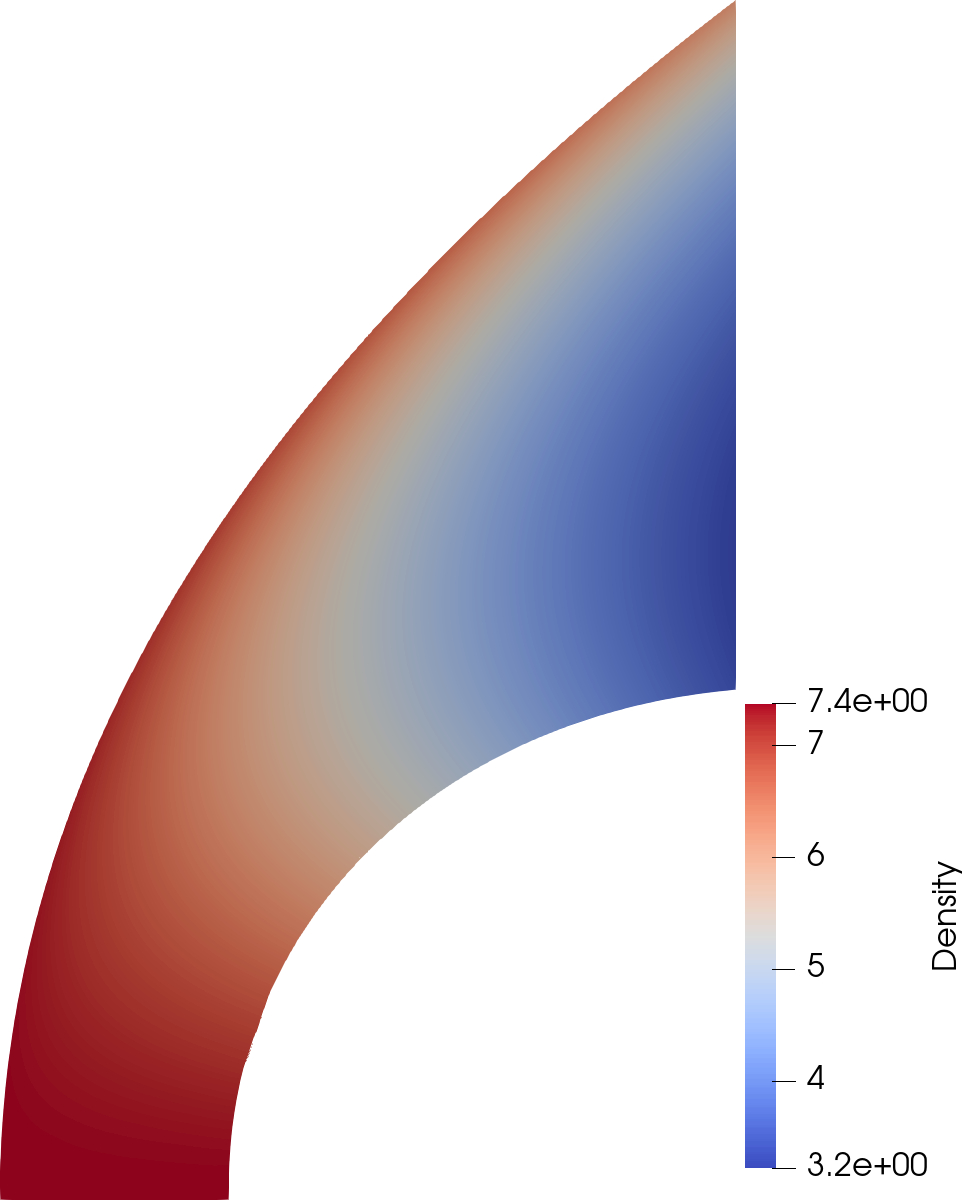}
\caption{Density variation for Mach number $6.0$ obtained using \eqref{finalP}, \eqref{finalY}, and \eqref{densityFormula} with shock shape given in \eqref{shapeBetterApprox} with parameters given in Table \ref{table:1} for 200 streamlines.}
\label{fig:M6p0Density}
\end{center}
\end{figure}

\begin{figure}[htbp]
\begin{center}
\includegraphics[scale=0.22]{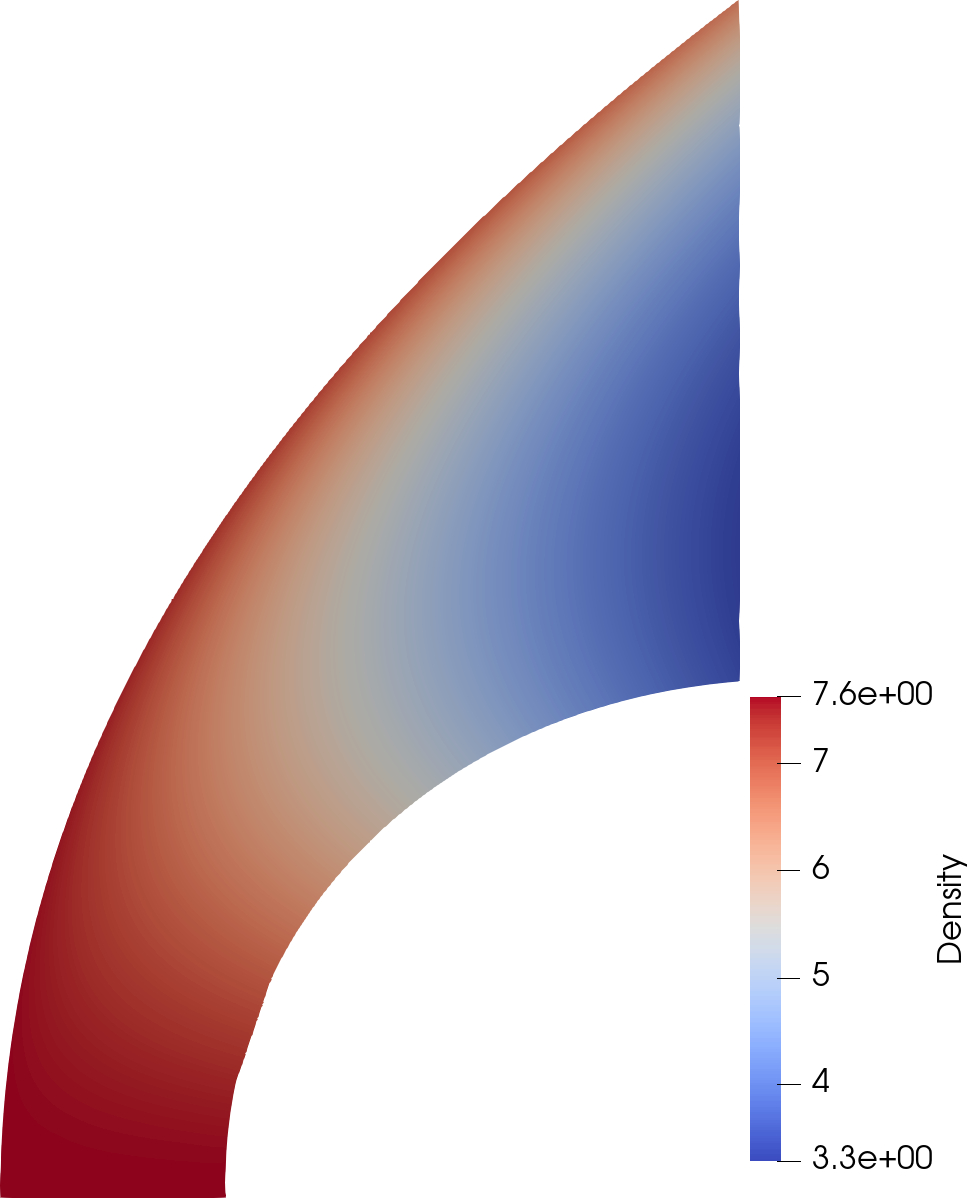}
\caption{Density variation for Mach number $7.0$ obtained using \eqref{finalP}, \eqref{finalY}, and \eqref{densityFormula} with shock shape given in \eqref{shapeBetterApprox} with parameters given in Table \ref{table:1} for 200 streamlines.}
\label{fig:M7p0Density}
\end{center}
\end{figure}

\begin{figure}[htbp]
\begin{center}
\includegraphics[scale=0.22]{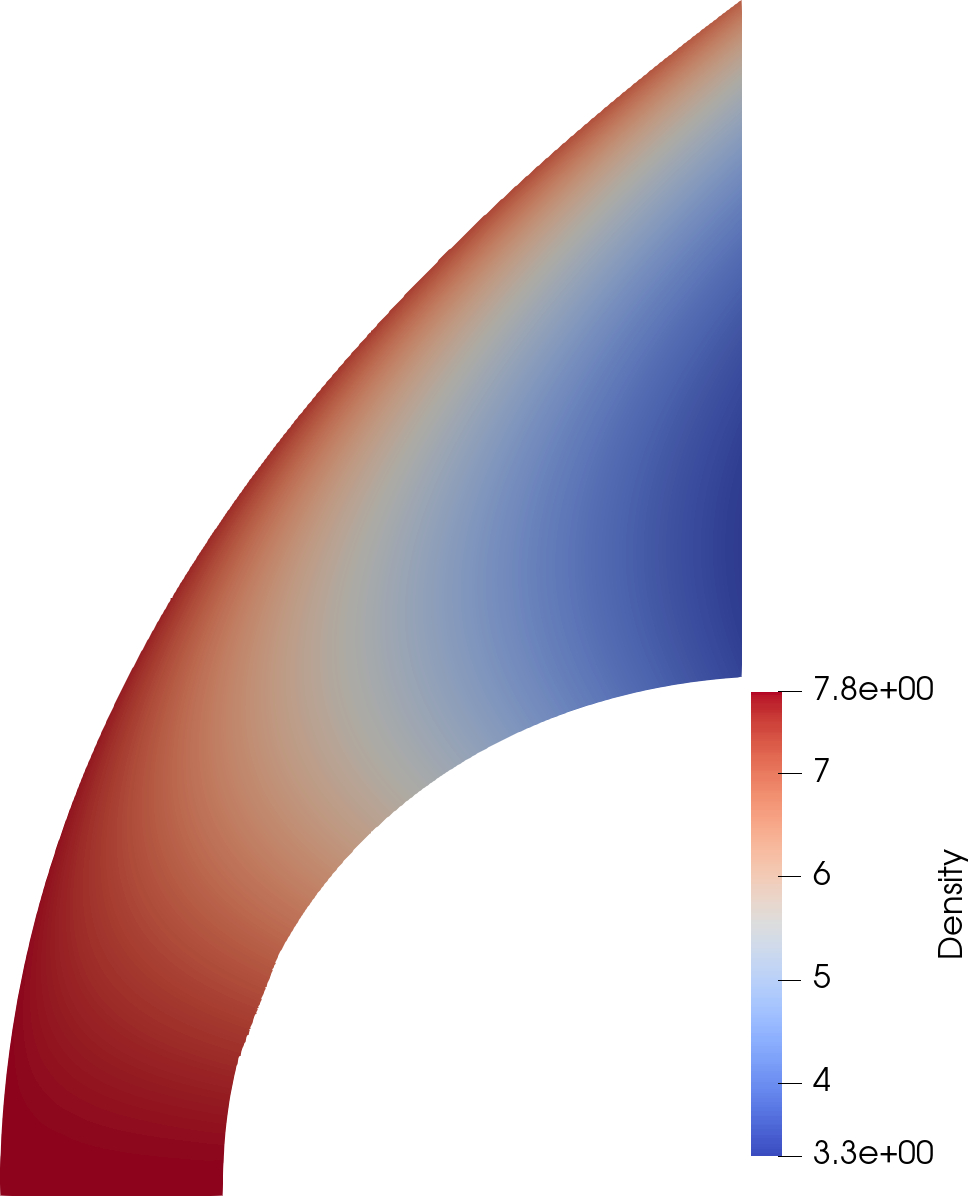}
\caption{Density variation for Mach number $8.0$ obtained using \eqref{finalP}, \eqref{finalY}, and \eqref{densityFormula} with shock shape given in \eqref{shapeBetterApprox} with parameters given in Table \ref{table:1} for 200 streamlines.}
\label{fig:M8p0Density}
\end{center}
\end{figure}

\noindent After using this shock shape, we now calculate the full flow field between the shock and the body on the streamlines by evaluating equations \eqref{finalP}, \eqref{finalY}, and \eqref{densityFormula}. We have used 200 streamlines for each Mach number. The density variations for Mach numbers $4.0$, $5.0$, $6.0$, $7.0$, and $8.0$, obtained with this shock shape are shown in Figures \ref{fig:M4p0Density}, \ref{fig:M5p0Density}, \ref{fig:M6p0Density}, \ref{fig:M7p0Density}, and \ref{fig:M8p0Density} respectively. All the plots are obtained by plotting the point data and using a delaunay triangulation to convert them into cells to cover the whole domain using paraview \cite{paraview}.
\\

\begin{figure}[htbp]
\begin{center}
\includegraphics[scale=0.22]{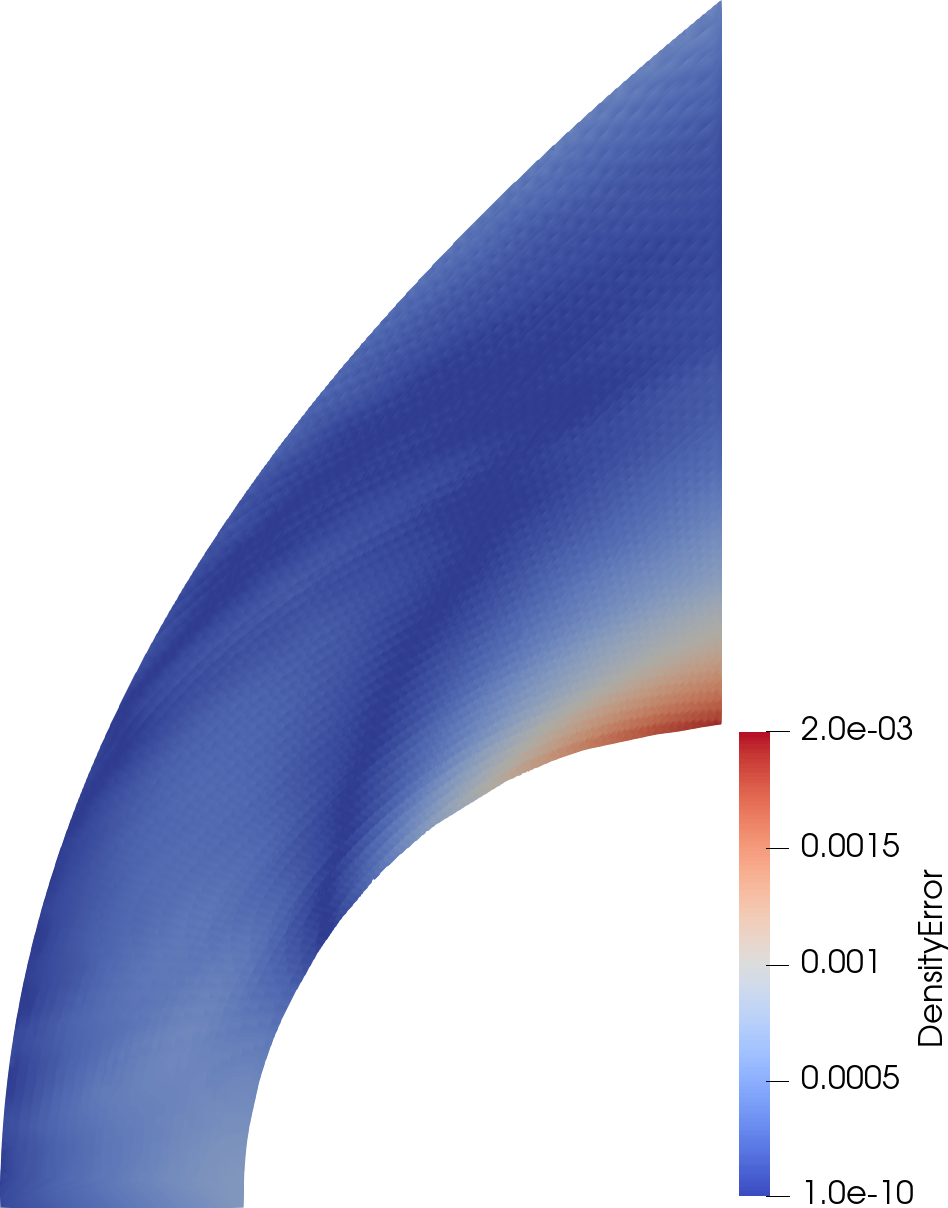}
\caption{Density error obtained for Mach number $5.0$ using a fifth order accurate (formally called $\mathbf{P}^{4}$ based) discontinuous Galerkin method along with overset grids using 9976 elements for accurate shock capturing \cite{srspkmr4}.}
\label{fig:M5p0DensityError}
\end{center}
\end{figure}

\begin{figure}[htbp]
\begin{center}
\includegraphics[scale=0.22]{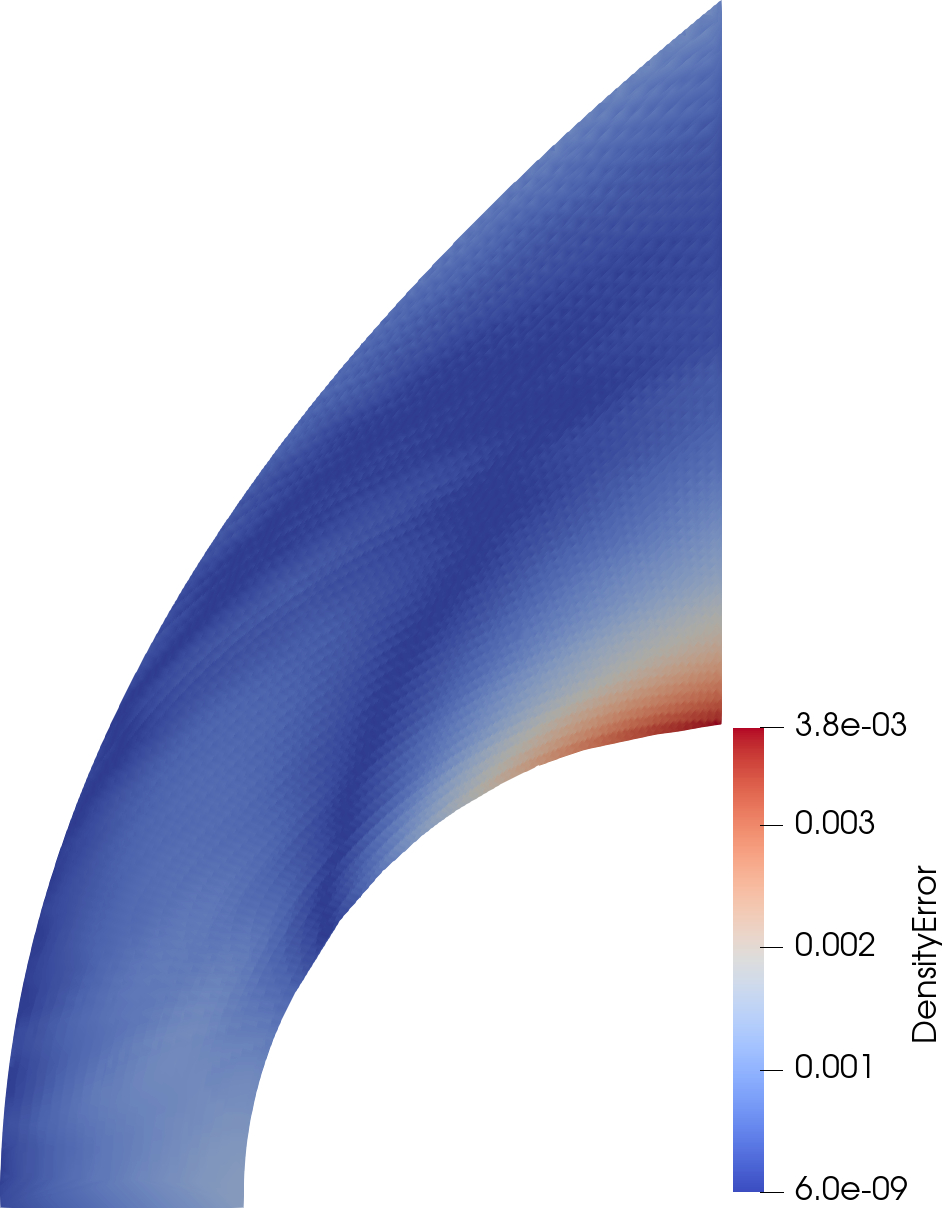}
\caption{Density error obtained for Mach number $6.0$ using a fifth order accurate (formally called $\mathbf{P}^{4}$ based) discontinuous Galerkin method along with overset grids using 9976 elements for accurate shock capturing \cite{srspkmr4}.}
\label{fig:M6p0DensityError}
\end{center}
\end{figure}

\begin{figure}[htbp]
\begin{center}
\includegraphics[scale=0.22]{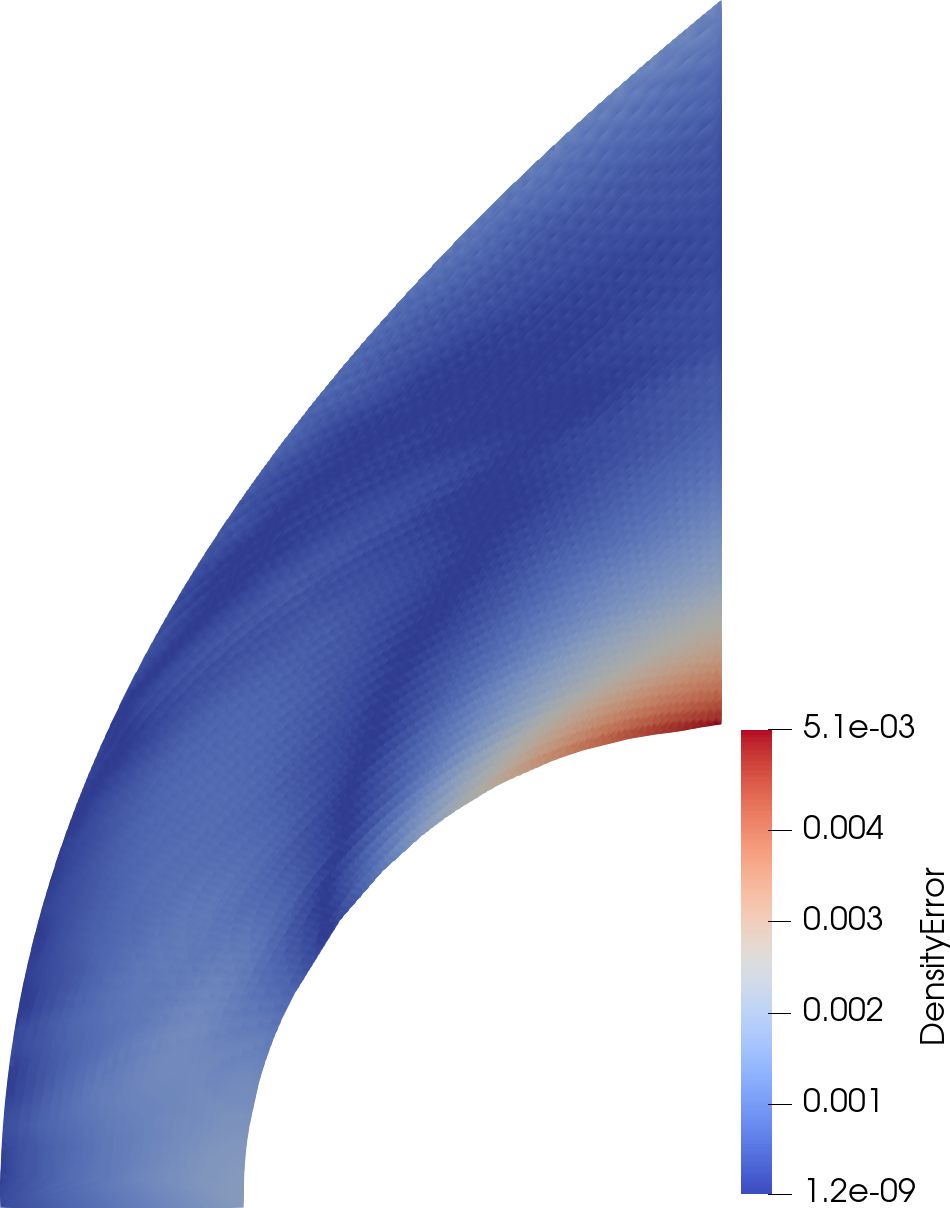}
\caption{Density error obtained for Mach number $7.0$ using a fifth order accurate (formally called $\mathbf{P}^{4}$ based) discontinuous Galerkin method along with overset grids using 9976 elements for accurate shock capturing \cite{srspkmr4}.}
\label{fig:M7p0DensityError}
\end{center}
\end{figure}

\begin{figure}[htbp]
\begin{center}
\includegraphics[scale=0.22]{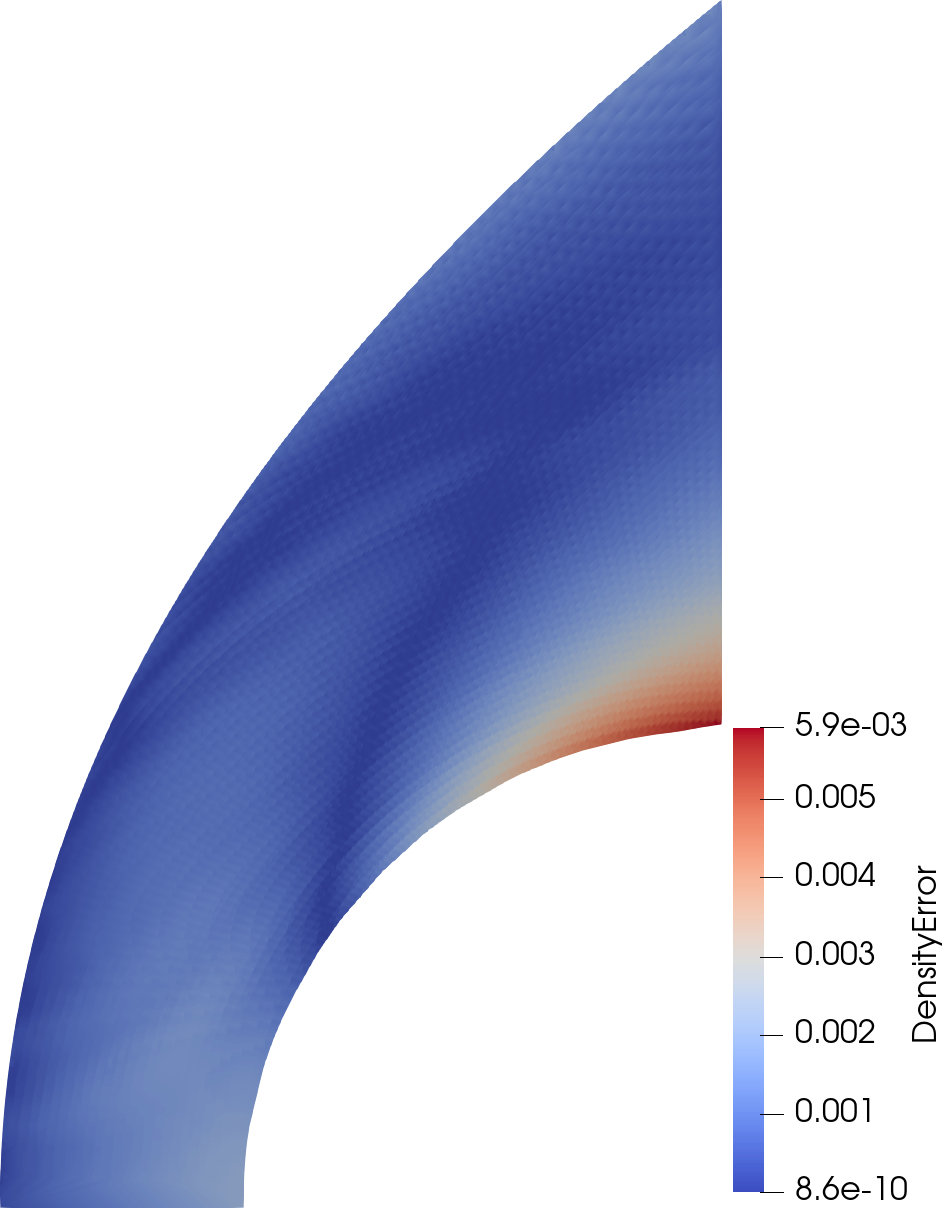}
\caption{Density error obtained for Mach number $8.0$ using a fifth order accurate (formally called $\mathbf{P}^{4}$ based) discontinuous Galerkin method along with overset grids using 9976 elements for accurate shock capturing \cite{srspkmr4}.}
\label{fig:M8p0DensityError}
\end{center}
\end{figure}

\noindent We now validate our procedure by using a numerically calculated solution. We have used a fifth order accurate (formally called $\mathbf{P}^{4}$ based) discontinuous Galerkin method (DGM) along with overset grids  with 9976 elements for accurate shock capturing \cite{srspkmr4} to solve the two-dimensional Euler equations over a circular cylinder and used this solution to calculate the density error in our solution. The formulation and validation of our solver using DGM, along with a brief outline of the shock capturing procedure are given in Appendix \ref{app:DGMFormulationValidation}. The density errors obtained using this fifth order accurate numerical solution for Mach number $4.0$, $5.0$, $6.0$, $7.0$, and $8.0$, are shown in Figures \ref{fig:M4p0Density}, \ref{fig:M5p0DensityError}, \ref{fig:M6p0DensityError}, \ref{fig:M7p0DensityError}, and \ref{fig:M8p0DensityError} respectively. We don't show the solution obtained by the discontinuous Galerkin method as it is quite similar to the solution obtained using our procedure and shown in Figures \ref{fig:M4p0Density} - \ref{fig:M8p0Density}.
\\
\begin{figure}[htbp]
\begin{center}
\includegraphics[scale=0.22]{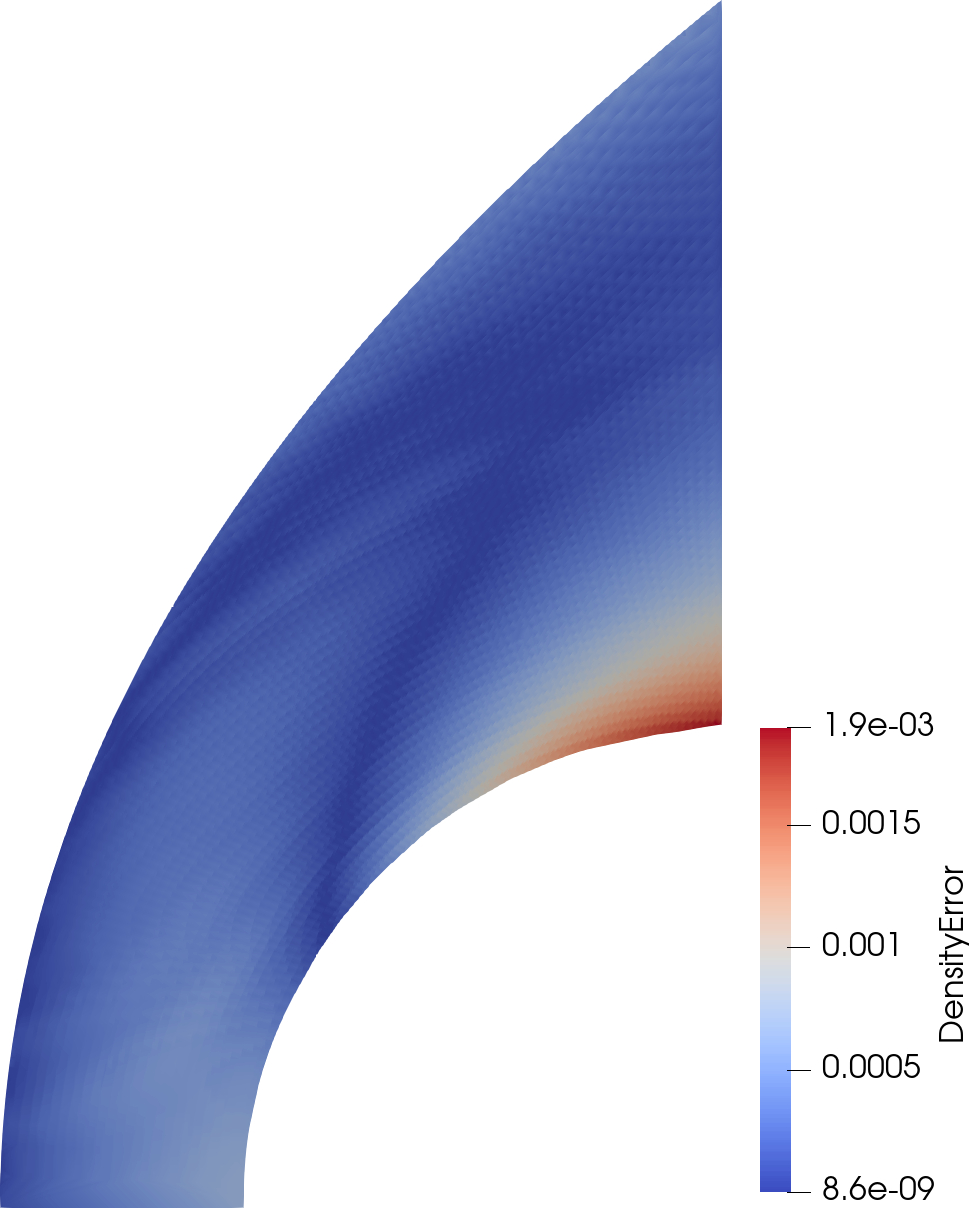}
\caption{Density error obtained for Mach number $4.0$ using a fifth order accurate (formally called $\mathbf{P}^{4}$ based) discontinuous Galerkin method along with overset grids using 9976 elements for accurate shock capturing \cite{srspkmr4}.}
\label{fig:M4p0DensityError}
\end{center}
\end{figure}

\noindent We now explain the choice of the selected Mach numbers. We have calculated the maximum density error obtained using the procedure in Section \ref{formulation} and tabulated them in Table \ref{table:2}. We also show the value of $\rho_{\infty}/\rho$ across a normal shock in Table \ref{table:2} for reference. We can clearly see that the error is increasing quite fast when we go below Mach $4.0$. In fact, it has increased by two orders of magnitude when we go from Mach $4.0$ to Mach $3.5$. The reason for this can be the first assumption for this procedure given by equation \eqref{assump1} where we have assumed that the density immediately behind a strong shock is much larger than in front of the shock. The values of $\rho_{\infty}/\rho$ (across a normal shock) given in Table \ref{table:2} decrease when we increase the Mach number and we probably have to cross a threshold for this procedure to work well which can be between Mach numbers 3.5 and 4.0 as seen from Table \ref{table:2}. For this reason, we have selected Mach numbers $\ge 4$ to illustrate our solution.
\\
\begin{table}[htbp]
\centering
\begin{tabular}{|c|c|c|}
\hline
Mach Number & Maximum Density Error & $\rho_{\infty}/\rho$ across a normal shock \\
\hline
3.5 & 1.1E-01 & 0.2347 \\
\hline
3.8 & 8.9E-02 & 0.2244 \\
\hline
3.9 & 2.1E-02 & 0.2214 \\
\hline
4.0 & 1.9E-03 & 0.2188 \\
\hline
5.0 & 2.0E-03 & 0.2 \\
\hline
6.0 & 3.8E-03 & 0.1898 \\
\hline
7.0 & 5.1E-03 & 0.1837 \\
\hline
8.0 & 5.9E-03 & 0.1797 \\
\hline
\end{tabular}
\caption{Maximum density error obtained using the procedure given in Section \ref{formulation} for various Mach numbers. For these Mach numbers, the density ratio ($\rho_{\infty}/\rho$) across a normal shock is also shown. The error is calculated using a fifth order accurate discontinuous Galerkin method along with overset grids using 9976 elements for accurate shock capturing \cite{srspkmr4}}
\label{table:2}
\end{table}

\noindent To further validate the solution obtained, we calculate the functional - integrated surface pressure on the body of the cylinder and compare it with the numerical solution obtained. The error in integrated surface pressure is tabulated in Table \ref{table:3}. From the values obtained, we can see that the error in the integrated surface pressure is smaller than the maximum density error. This shows that the procedure is quite accurate near the body (except at some locations as will be discussed later) in predicting the surface pressure which is an important quantity.
\\
\begin{table}[htbp]
\centering
\begin{tabular}{|c|c|}
\hline
Mach Number & Error in integrated surface pressure \\
\hline
4.0 & 3.45E-04 \\
\hline
5.0 & 2.75E-04 \\
\hline
6.0 & 2.15E-04  \\
\hline
7.0 & 2.03E-04 \\
\hline
8.0 & 1.88E-04 \\
\hline
\end{tabular}
\caption{Error in the functional - integrated surface pressure for various Mach numbers using the shock shape given by $af+bf^{2}$. The error is calculated using a fifth order accurate discontinuous Galerkin method along with overset grids using 9976 elements for accurate shock capturing \cite{srspkmr4}}
\label{table:3}
\end{table}

\begin{figure}[htbp]
\begin{center}
\includegraphics[scale=0.22]{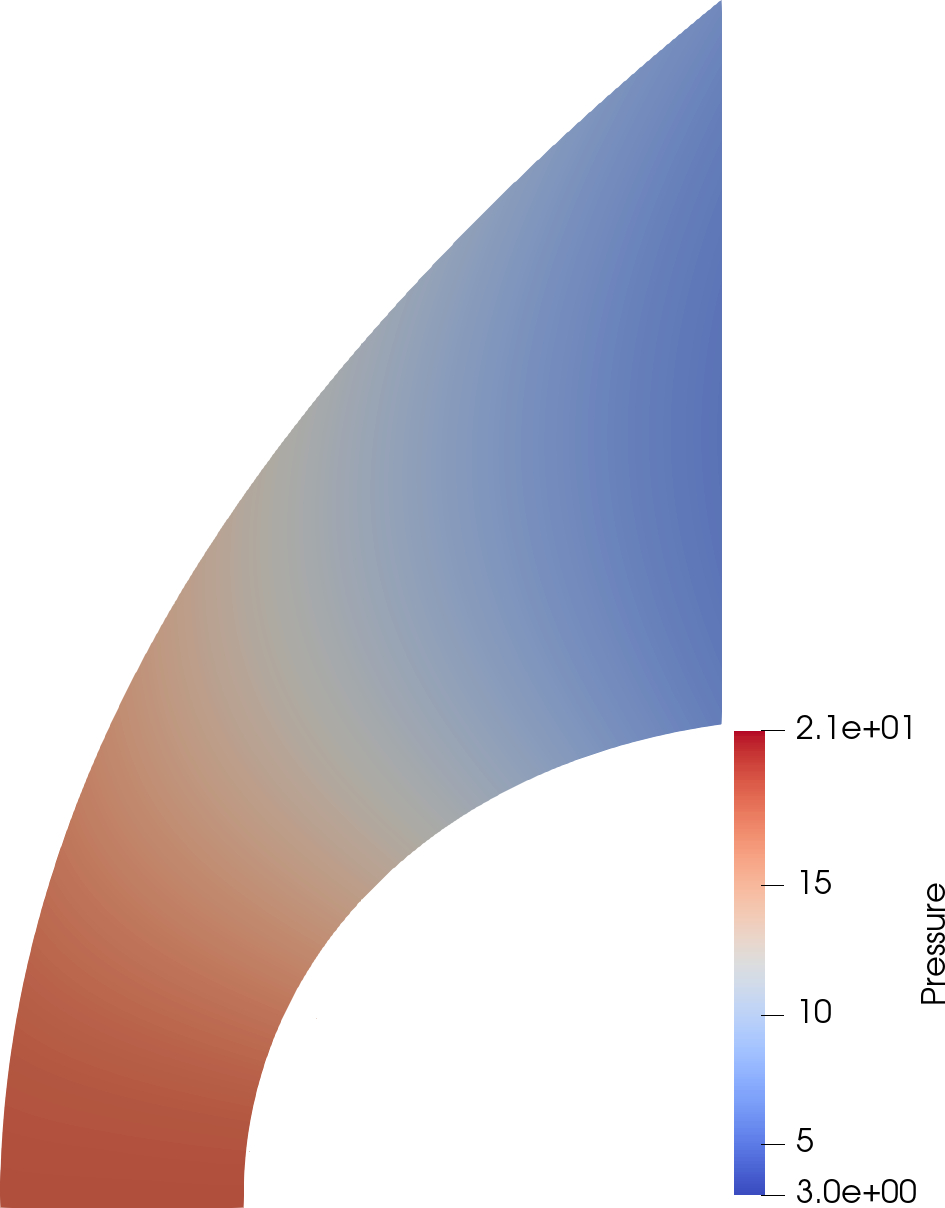}
\caption{Pressure variation for Mach number $4.0$ with shock shape given in \eqref{shapeBetterApprox} with parameters given in Table \ref{table:1} for 200 streamlines.}
\label{fig:M4p0Pressure}
\end{center}
\end{figure}

\noindent \textbf{Discussion of results:} We can see that the results obtained for each Mach number are very accurate and the maximum error in comparison with a numerical solution for each case is only about the order of $10^{-3}$ (using a shock shape given by $af+bf^2$). For all the Mach numbers, the error pattern is similar and the maximum error occurs near the upper portion of the body. To explain this, we show the pressure solution for Mach $4.0$ in Figure \ref{fig:M4p0Pressure}. We also remember the assumption that the pressure at the point Q (Figure \ref{fig:ShockBody}) of the disturbed flow field is not much smaller than the pressure at point N (see Figure \ref{fig:ShockBody}). This assumption is given by equation \eqref{assump2}. From Figure \ref{fig:M4p0Pressure}, we can see that for the points that are near the upper portion of the body, the pressure ratio $\hat{P}/P$ is larger compared to the other points. Though we showed the pressure variation for Mach 4.0 only, this happens to be true for Mach numbers 5.0, 6.0, 7.0, and 8.0 as well. This can be why the error obtained is the largest in that region as we are slightly violating an assumption with which we derived the solution.
\\
\\

\begin{figure}[htbp]
\begin{center}
\includegraphics[scale=1.0]{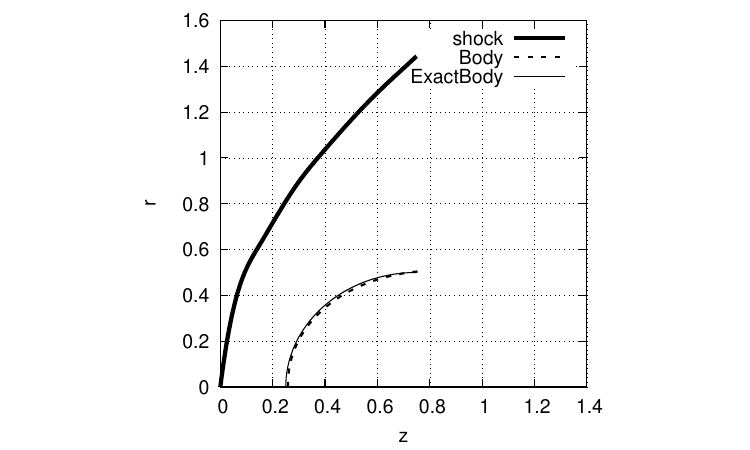}
\caption{The shock using the cubic spline from the numerical solution (thick solid line), the body obtained from \eqref{YatZero} (for plane flow) (thin dashed line), and the actual shape of the body (thin solid line) for Mach number $M=4.0$. Radius of the circular cylinder is $0.5$}
\label{fig:M4p0ShockBodyShapeNewNumeric}
\end{center}
\end{figure}

\begin{figure}[htbp]
\begin{center}
\includegraphics[scale=0.22]{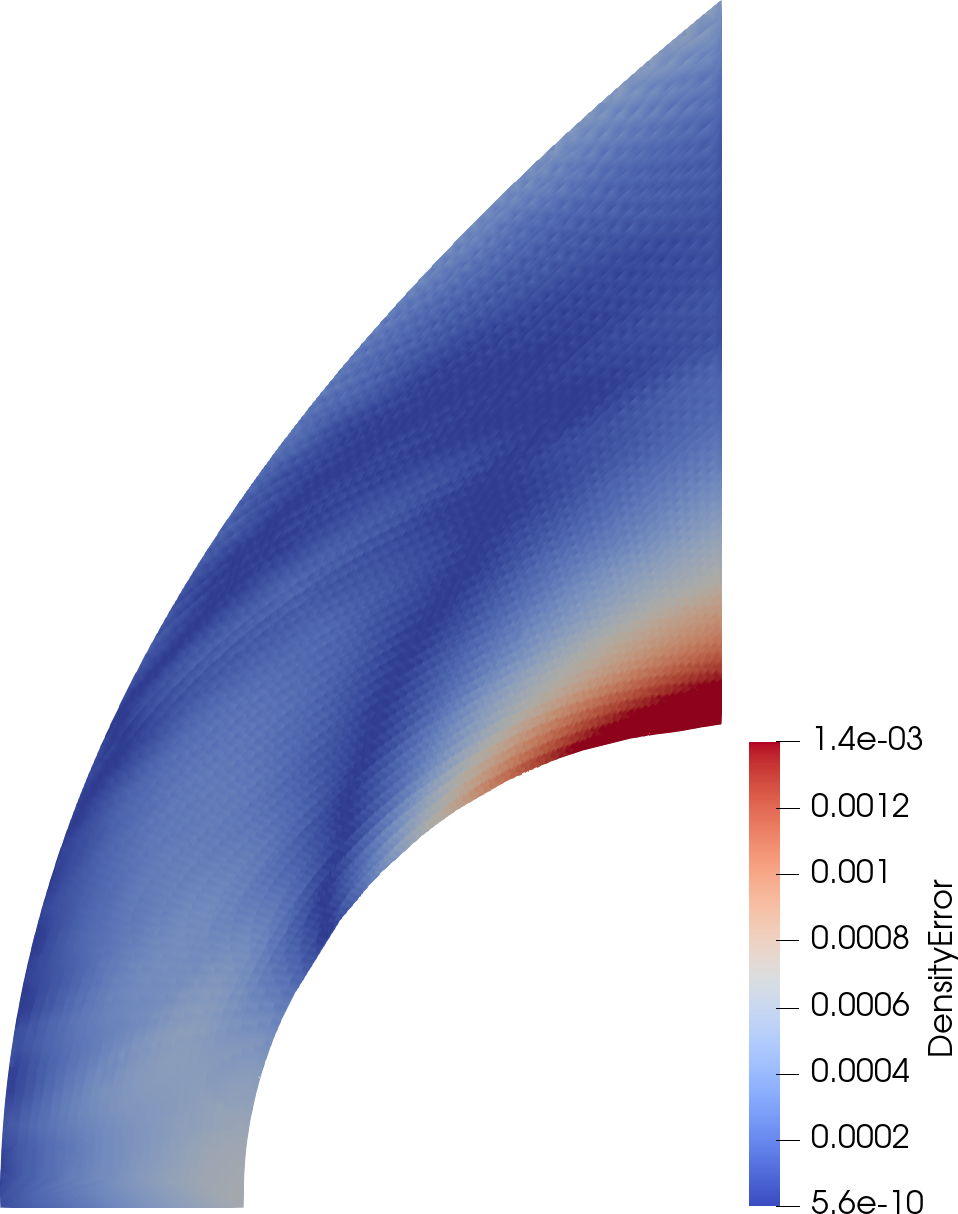}
\caption{Density error obtained for Mach number $4.0$ using the cubic spline shock shape from the numerical solution using a fifth order accurate (formally called $\mathbf{P}^{4}$ based) discontinuous Galerkin method along with overset grids using 9976 elements for accurate shock capturing \cite{srspkmr4}.}
\label{fig:M4p0DensityErrorNumeric}
\end{center}
\end{figure}

\noindent We have also changed the shock shape further by looking at equations of the form $af+bf^{2}+cf^{3}, af+bf^{2}+cf^{3}+df^{4}$ with $f$ defined by \eqref{fDef}. We have optimized the parameters $z_{0}$, $a$, $b$, $c$, and $d$ as mentioned in Section \ref{formulation}. We tabulate the maximum density errors for shock shapes with polynomial expansions in $f$ till fourth degree in Table \ref{table:3}. From this table, we can see that the solution accuracy increases only till the second degree approximation in $f$ and doesn't improve much for higher degree approximations.
\\
\\
\noindent As a final shock shape, we consider the accurate numerical solution obtained in \cite{srspkmr4} where the shock has been captured aligned to a grid line. We use that solution and obtain a cubic spline for the shock shape. Using the cubic spline as the shock shape, we obtain the body and the full solution in the shock layer using the method given in Section \ref{formulation}. The shape of the body obtained along with the shock and the exact shape of the body are shown in Figure \ref{fig:M4p0ShockBodyShapeNewNumeric} for Mach $4.0$ flow. We can see that the shape of the body obtained in this fashion is quite accurate. The density error for this solution is shown in Figure \ref{fig:M4p0DensityErrorNumeric}. We can see that the maximum density error is about $1.4\times 10^{-3}$. We also tabulate the maximum density error obtained in this fashion for various Mach numbers in Table \ref{table:4}.
\\

\begin{table}[htbp]
\centering
\begin{tabular}{|c|c|c|c|c|c|}
\hline
 & \multicolumn{5}{|c|}{Maximum density error with shock shape using Schneider's method for plane flow:} \\
\cline{2-6}
\makecell{Mach \\ number} & $f$ & \makecell{$af+bf^{2}$ \\ (shock shape \\ used for all \\ our calculations)}  & $af+bf^{2}+cf^{3}$ & $af+bf^{2}+cf^{3}+df^{4}$ & \makecell{cubic spline \\ obtained from \\ numerical solution}\\
 \hline
4.0 & 1.8E-02 & 1.9E-03 & 1.8E-03 & 1.8E-03 & 1.4E-03 \\
\hline
5.0 & 1.9E-02 & 2.0E-03 & 1.7E-03 & 1.7E-03 & 1.1E-03 \\
\hline
6.0 & 2.1E-02 & 3.8E-03 & 3.5E-03 & 3.4E-03 & 2.4E-03 \\
\hline
7.0 & 2.2E-02 & 5.1E-03 & 4.9E-03 & 4.8E-03 & 3.6E-03 \\
\hline
8.0 & 2.8E-02 & 5.9E-03 & 5.8E-03 & 5.5E-03 & 4.3E-03 \\
\hline
\end{tabular}
\caption{Plane flow: Maximum density error obtained for various Mach numbers with various shock shapes obtained using an expansion in terms of $f$ and also using a cubic spline obtained from the numerical solution. The error is calculated using a fifth order accurate discontinuous Galerkin method along with overset grids using 9976 elements for accurate shock capturing \cite{srspkmr4}}
\label{table:4}
\end{table}

\noindent From these results, we can conclude that using this kind of a shock shape, this might be the best accuracy we can get using this analytical method with the stated approximations and assumptions.
\\
\\
\noindent To see how well the shock shapes assumed work in a typical shock fitting solver, we use the shock shapes $f$, $af+bf^{2}$, $af+bf^{2}+cf^{3}$, $af+bf^{2}+cf^{3}+df^{4}$ and optimise the various parameters ($z_{0}, a, b, c, d$ as required) based on the boundary condition at the circular cylinder (the no-penetration condition). We get different values of the optimised parameters to those obtained using Schneider's method (as expected) and we tabulate the maximum density error for each of the shock shapes in Table \ref{table:5}. Here unlike the case for Schneider's method, we see that the maximum density error decreases reasonably as we expand the shock shape in series expansion. The decrease is quite prominent (two orders of magnitude) when we go from $f$ to $af+bf^{2}$. The final column of the table is given just to show that the shock shape obtained from numerical solution as a cubic spline solves the problem very accurately. All these results show that the choice of shock shape used for our solution is quite good.

\begin{table}[htbp]
\centering
\begin{tabular}{|c|c|c|c|c|c|}
\hline
 & \multicolumn{5}{|c|}{Maximum density error with shock shape with shock fitting using DGM for plane flow:} \\
\cline{2-6}
\makecell{Mach \\ number} & $f$ & $af+bf^{2}$  & $af+bf^{2}+cf^{3}$ & $af+bf^{2}+cf^{3}+df^{4}$ & \makecell{cubic spline \\ obtained from \\ numerical solution}\\
 \hline
4.0 & 3.3E-06 & 5.5E-08 & 3.7E-09 & 4.1E-10 & 1.8E-16 \\
\hline
5.0 & 2.4E-06 & 2.3E-08 & 4.4E-09 & 3.6E-10 & 2.1E-16 \\
\hline
6.0 & 4.1E-06 & 4.2E-08 & 3.2E-09 & 3.9E-10 & 1.3E-16 \\
\hline
7.0 & 5.2E-06 & 6.3E-08 & 8.3E-09 & 6.7E-10 & 2.2E-16 \\
\hline
8.0 & 4.8E-06 & 7.7E-08 & 6.9E-09 & 7.8E-10 & 3.2E-16 \\
\hline
\end{tabular}
\caption{Plane flow: Maximum density error obtained for various Mach numbers using shock fitting with DGM with various shock shapes obtained using an expansion in terms of $f$ and also using a cubic spline obtained from the numerical solution. The error is calculated using a fifth order accurate discontinuous Galerkin method along with overset grids using 9976 elements for accurate shock capturing \cite{srspkmr4}}
\label{table:5}
\end{table}

\subsection{Axisymmetric Flow:}\label{axisymmetric}

\noindent We have also tested Schneider's method given in Section \ref{formulation} for the axisymmetric case. To obtain an accurate numerical solution for this case, we have converted our plane flow solver into an axisymmetric flow solver using the method given by Yu \cite{stYu}. This gives us a fifth order accurate discontinuous Galerkin solution for axisymmetric flow. We have used this numerical solution with 7500 elements to calculate the errors in the analytical solution. We used the same shock shape given by Moeckel with 

\begin{equation}\label{axisymmShockShape}
 r = f
\end{equation}

\noindent where $f$ is given by equation \eqref{fDef}. Here, $r$ is the radial coordinate and $z$ is the longitudinal coordinate.
\\
\\
\noindent We now use the shock shapes $f$, $af+bf^{2}$, $af+bf^{2}+cf^{3}$, $af+bf^{2}+cf^{3}+df^{4}$ and optimise the various parameters ($z_{0}, a, b, c, d$ as required) again using the error between the exact body shape and the body shape obtained using the method given in Section \ref{formulation}. We get different values of the optimised parameters to those obtained for plane flow (as expected) and we tabulate the maximum density error for each of the shock shapes in Table \ref{table:6}. Here, we can see that the pattern for maximum density error is similar to that of plane flow decreasing from shock shapes $f$ to $af+bf^{2}$ and remaining almost the same for further series expansions. We also show the error obtained using a shock shock obtained from numerical solution in the last column of Table \ref{table:6}. From these results, we can again conclude that using this kind of a shock shape, this might be the best accuracy we can get for axisymmetric flow using the analytical method given in Section \ref{formulation}.

\begin{table}[htbp]
\centering
\begin{tabular}{|c|c|c|c|c|c|}
\hline
 & \multicolumn{5}{|c|}{Maximum density error with shock shape using Schneider's method for axisymmetric flow:} \\
\cline{2-6}
\makecell{Mach \\ number} & $f$ & $af+bf^{2}$ & $af+bf^{2}+cf^{3}$ & $af+bf^{2}+cf^{3}+df^{4}$ & \makecell{cubic spline \\ obtained from \\ numerical solution}\\
 \hline
4.0 & 8.5E-02 & 4.3E-03 & 4.2E-03 & 3.8E-03 & 1.6E-03 \\
\hline
5.0 & 8.3E-02 & 3.4E-03 & 3.1E-03 & 2.8E-03 & 1.3E-03 \\
\hline
6.0 & 7.4E-02 & 3.2E-03 & 2.8E-03 & 2.6E-03 & 1.1E-03 \\
\hline
7.0 & 7.3E-02 & 2.8E-03 & 2.5E-03 & 2.3E-03 & 1.2E-03 \\
\hline
8.0 & 6.6E-02 & 2.6E-03 & 2.3E-03 & 2.2E-03 & 1.7E-03 \\
\hline
\end{tabular}
\caption{Axisymmetric Flow: Maximum density error obtained for various Mach numbers with various shock shapes obtained using an expansion in terms of $f$ and also using a cubic spline obtained from the numerical solution. The error is calculated using a fifth order accurate discontinuous Galerkin solver for axisymmetric flow with 7500 elements.}
\label{table:6}
\end{table}

\section{Conclusion}\label{conc}

We have calculated the high supersonic flow over a circular cylinder using Schneider's inverse method \cite{schneider}. We assumed the shock shape of a hyperbola in the form given by Moeckel in \cite{moeckel} and optimized the shock shape by minimising the error between the shape of the body obtained using Schneider's method and the actual shape of the body. We have improved this shock shape further by writing it as a series expansion ($af+bf^{2}$, where $f$ gives the shock shape of Moeckel) with unknown coefficients. We have determined these coefficients using the same optimization and obtained a much better approximation for the shock shape. We have used this shock shape to find the solution of the flow in the shock layer using Schneider's method by integrating the equations of motion using the stream function. We have compared the solution obtained with a numerical solution calculated using a fifth order accurate discontinuous Galerkin method which captures the shock very accurately using overset grids \cite{srspkmr4}. We have found that the maximum error in density is only of the order of $10^{-3}$ (with the shock shape $af+bf^{2}$) demonstrating the accuracy of the solution method. We also found that using more terms in the series expansion for the shock shape (like $af+bf^{2}+cf^{3}$, $af+bf^{2}+cf^{3}+df^{4}$) does not improve the accuracy of the solution. We also used a cubic spline obtained from the numerical solution to find the solution and the maximum error in density is again of the order of $10^{-3}$. This suggests that this might be the best accuracy we can get using this analytical method. We have also used the shock shape given by Moeckel and used it in a series expansion for axisymmetric flow case and optimised the parameters. The density errors again follow the same pattern of plane flow and the maximum density error is again of the order of $10^{-3}$. 
\begin{appendices}
\section{Formulation and Validation of discontinuous Galerkin Method}\label{app:DGMFormulationValidation}
 
\noindent Consider the Euler equations in conservative form as given by

\begin{equation}\label{2dEulerEquations}
\frac{\partial \textbf{Q}}{\partial t} + \frac{\partial \textbf{F(Q)}}{\partial x} + \frac{\partial \textbf{G(Q)}}{\partial y} = 0 \quad \text{in the domain} \quad \Omega
\end{equation}
\noindent where $\textbf{Q} = (\rho, \rho u, \rho v, E)^{T}$, $\textbf{F(Q)}=u\textbf{Q} + (0, P, 0, Pu)^{T}$ and $\textbf{G(Q)}=v\textbf{Q} + (0, 0, P, Pv)^{T}$ with $P = (\gamma -1)(E-\frac{1}{2}\rho (u^{2}+v^{2}))$ and $\gamma = 1.4$. Here, $\rho$ is the density, $(u,v)$ is the velocity, $E$ is the total energy and $P$ is the pressure. We approximate the domain $\Omega$ by $K$ non overlapping elements given by $\Omega_{k}$. 
\\
\\
We look at solving \eqref{2dEulerEquations} using the discontinuous Galerkin method (DGM). We approximate the local solution in an element $\Omega_{k}$, where $k$ is the element number, as a polynomial of order $N$ which is given by:

\begin{equation}\label{modalApprox}
 Q_{h}^{k}(r,s) = \sum_{i=0}^{N_{p}-1} Q_{i}^{k} \psi_{i}(r,s)
\end{equation}

\noindent where $N_{p}=(N+1)(N+1)$ and $r$ and $s$ are the local coordinates. Here, the subscript $i$ represents the particular degree of freedom, $h$ represents the grid size, and the superscript $k$ is the element number. The polynomial basis used ($\psi_{i}(r,s)$) is the tensor product orthonormalized Legendre polynomials of degree $N$. The number of degrees of freedom are given by $N_{p}=(N+1)(N+1)$. Now, using $\psi_{j}(r,s)$ as the test function, the weak form of the equation \eqref{2dEulerEquations} is obtained as

\begin{equation}\label{weakFormScheme}
 \sum_{i=0}^{N_{p}-1} \frac{\partial Q_{i}^{k}}{\partial t} \int_{\Omega_{k}} \psi_{i} \psi_{j} d\Omega + \int_{\partial \Omega_{k}} \hat{F} \psi_{j} ds - \int_{\Omega_{k}} \vec{F}.\nabla \psi_{j} d\Omega = 0 \quad j = 0,\ldots,N_{p}-1
\end{equation}

\noindent where $\partial \Omega_{k}$ is the boundary of $\Omega_{k}$, $\vec{F} = (\textbf{F(Q)},\textbf{G(Q)})$ and $\hat{F} = \bar{F^{*}}.\hat{n}$ where $\bar{F^{*}}$ is the monotone numerical flux at the interface which is calculated using an exact or approximate Riemann solver and $\hat{n}$ is the unit outward normal. We have used the Lax-Friedrichs numerical flux for all our calculations unless otherwise specified. This is termed to be $\mathbf{P}^{N}$ based discontinuous Galerkin method.
\\
\\
\noindent Equation \eqref{weakFormScheme} is integrated using an appropriate Gauss Legendre quadrature and is discretized in time by using the fifth order Runge-Kutta time discretization given in \cite{butcher} unless otherwise specified. To control spurious oscillations which occur near discontinuities, a limiter is used with a troubled cell indicator (a shock detector). We have used the KXRCF troubled cell indicator \cite{kxrcf} and the compact subcell WENO (CSWENO) limiter proposed in \cite{srspkmr1} for all our calculations. To capture the shock accurately, we have used overset grids as proposed in \cite{srspkmr4}. The communication between overset grids occurs using the procedure given in \cite{srspkmr3}. The shock capturing procedure with overset grids is briefly given below.
\\
\\
\noindent \textbf{Shock capturing procedure with overset grids:} The step by step procedure is as follows:
\\
\\
\noindent \textbf{Step 1:} Run the solver on a coarse grid with a given troubled cell indicator and limiter to steady state and obtain the solution.
\\
\\
\noindent \textbf{Step 2:} Look at the troubled cells (with a reliable shock detector \cite{kxrcf}) to locate the discontinuities (shocks) that occur in the solution.The troubled cells give us a good idea of the location of the discontinuities.
\\
\\
\noindent \textbf{Step 3:} Construct an overset grid conforming to the computational domain which is refined in a direction perpendicular to the discontinuities such that they are approximately parallel to a grid line.
\\
\\
\noindent \textbf{Step 4:} Using this overset grid, we rerun the solver with the coarse grid solution as the initial condition. While running the solver, we use the troubled cell indicator and the limiter only on the overset grid. We also use a high resolution numerical flux on the overset grid to capture the shock accurately. We have used the SLAU2 \cite{ks3} numerical flux in the overset grid and the less expensive Lax-Friedrichs flux elsewhere. Using this procedure, we obtain a more accurate solution with the discontinuities approximately aligned to a grid line.
\\
\\
\noindent \textbf{Validation of the solver using Isentropic Vortex Problem:} Consider the two-dimensional Euler equations given by equation \eqref{2dEulerEquations} in the domain $[0,10] \times [-5,5]$ for the Isentropic Euler Vortex problem suggested in \cite{shu1} as a test case. The analytical solution is given by: \\ $\rho = \left(1 -  \left(\frac{\gamma - 1}{16\gamma \pi^{2}}\right)\beta^{2} e^{2(1-r^{2})}\right)^{\frac{1}{\gamma-1}}$, $u = 1 - \beta e^{(1-r^{2})} \frac{y-y_{0}}{2\pi}$, $v = \beta e^{(1-r^{2})} \frac{x-x_{0}-t}{2\pi}$, and $p = \rho^{\gamma}$, where $r=\sqrt{(x-x_{0}-t)^{2}+(y-y_{0})^{2}}$, $x_{0}=5$, $y_{0}=0$, $\beta=5$ and $\gamma = 1.4$. Initialization is done with the exact solution at $t=0$ and periodic boundary conditions are used at the edges of the domain. The solution is obtained for various orders making sure that the spatial (DGM) and temporal (Runge-Kutta) schemes are of the same order. The errors in density and numerical orders of accuracy are calculated at $t=2$ for various grid sizes and are presented in Table \ref{2DEulerEqnIsenVortexErrorsRec}. We can clearly see that the theoretical order of accuracy of the scheme is obtained by the DG solver.

\begin{table}[h]
\resizebox{1.0\textwidth}{!}{\begin{minipage}{1.0\textwidth}
\centering
\begin{tabular}{|c|c|c|c|c|c|c|c|c|c|}
\hline
 & $K\times K$ & $L_{2}$ error & Order & $L_{1}$ error & Order & $L_{\infty}$ error & Order \\ \hline
 \multirow{5}{*}{$P^{1}$} & 20$\times$20 & 1.58E-03 &  & 2.41E-03 &  & 1.77E-01 &  \\ \cline{2-8}
   & 40$\times$40 & 3.66E-04 & 2.11 & 5.54E-04 & 2.12 & 4.39E-02 & 2.01  \\ \cline{2-8}
   & 80$\times$80 & 7.28E-05 & 2.33 &  1.11E-04 & 2.32 & 1.11E-02 & 1.99 \\ \cline{2-8}
   & 160$\times$160 & 1.39E-05 & 2.39 & 2.15E-05 & 2.37 & 2.80E-03 & 1.98 \\ \hline
  \multirow{5}{*}{$P^{2}$} & 20$\times$20 & 5.19E-04 &  & 6.33E-04 &  & 9.42E-02 &   \\ \cline{2-8}
  & 40$\times$40 & 6.27E-05 & 3.05 & 7.70E-05 & 3.04 & 1.15E-02 & 3.03 \\ \cline{2-8}
  & 80$\times$80 & 7.11E-06 & 3.14 & 8.67E-06 & 3.15 & 1.30E-03 & 3.15 \\ \cline{2-8}
  & 160$\times$160 & 8.35E-07 & 3.09 & 1.01E-06 & 3.10 & 1.62E-04 & 3.00 \\ \hline
  \multirow{5}{*}{$P^{3}$} & 20$\times$20 & 3.30E-06 &  & 4.66E-06 &  & 8.21E-05 &  \\ \cline{2-8}
  & 40$\times$40 & 2.06E-07 & 4.00 & 2.93E-07 & 3.99 & 5.17E-06 & 3.99 \\ \cline{2-8}
  & 80$\times$80 & 1.27E-08 & 4.02 &  1.85E-08 & 3.99 & 3.30E-07 & 3.97   \\ \cline{2-8}
  & 160$\times$160 & 7.89E-10 & 4.01 & 1.12E-09 & 4.04 & 2.09E-08 & 3.98 \\ \hline
  \multirow{5}{*}{$P^{4}$} & 20$\times$20 & 1.76E-07 &  & 3.16E-07 &  & 9.83E-06 &  \\ \cline{2-8}
  & 40$\times$40 & 5.54E-09 & 4.99 & 1.01E-08 & 4.97 & 3.32E-07 & 4.89 \\ \cline{2-8}
  & 80$\times$80 & 1.78E-10 & 4.96 & 3.24E-10 & 4.96 & 1.10E-08 & 4.92 \\ \cline{2-8}
  & 160$\times$160 & 5.68E-12 & 4.97 & 1.05E-11 & 4.95 & 3.57E-10 & 4.94 \\ \hline
\end{tabular}
\end{minipage}}
\caption{2D Euler equations for the Isentropic Vortex problem with periodic boundary conditions, $t=2$, uniform mesh with $K\times K$ elements, $L_{1}$, $L_{2}$ and $L_{\infty}$ errors for $P^{1}$, $P^{2}$, $P^{3}$ and $P^{4}$ based DG}
\label{2DEulerEqnIsenVortexErrorsRec}
\end{table}
\end{appendices}
%
\bibliographystyle{ieeetr}
\bibliography{references}

\begin{thebibliography}{10}

\bibitem{schneider}
W.~Schneider, ``{A uniformly valid solution for the hypersonic flow past
  blunted bodies.},'' {\em Journal of Fluid Mechanics}, vol.~31, pp.~397--415,
  1968.

\bibitem{moeckel}
W.~Moeckel, ``{Approximate Method for Predicting Form and Location of Detached
  Shock Waves Ahead of Plane or Axially Symmetric Bodies.},'' {\em NACA
  Technical Note, No.1921}, 1949.

\bibitem{chester1}
W.~Chester, ``{Supersonic flow past a bluff body with a detached shock. Part I,
  Two-dimensional body.},'' {\em Journal of Fluid Mechanics}, vol.~1,
  pp.~353--365, 1956.

\bibitem{chester2}
W.~Chester, ``{Supersonic flow past a bluff body with a detached shock. Part
  II, Axisymmetrical body.},'' {\em Journal of Fluid Mechanics}, vol.~1,
  pp.~490--496, 1956.

\bibitem{freeman1}
N.~Freeman, ``{On the theory of hypersonic flow past plane and axially
  symmetric bluff bodies.},'' {\em Journal of Fluid Mechanics}, vol.~1,
  pp.~366--387, 1956.

\bibitem{hprobstein}
W.~Hayes and R.~Probstein, {\em {Hypersonic flow theory}}.
\newblock {Academic Press, New York}, 1966.

\bibitem{rasmussen1}
M.~Rasmussen, {\em {Hypersonic flow}}.
\newblock {Wiley}, 1994.

\bibitem{salas1}
M.~D. Salas, {\em {A Shock-Fitting Primer}}.
\newblock {CRC Applied Methematics and Nonlinear Science, Chapman and Hall, 1st
  edition}, 2009.

\bibitem{chou1}
Y.~Chou, ``{Radiative Coupled Viscous Flow with Massive Blowing}.,'' {\em NASA
  CR-2236}, 1973.

\bibitem{maslen1}
S.~Maslen, ``{Inviscid Flow about Blunted Cones of Large Opening Angle at Angle
  of Attack}.,'' {\em NASA CR-132652}, 1975.

\bibitem{nagaraja1}
K.~Nagaraja, ``{On the Shock Stand-off Distance in an Inviscid Hypersonic
  Source Flow Past Two-Dimensional Bluff Bodies}.,'' {\em ARL 71-0303}, 1971.

\bibitem{sb1}
M.~Schulreich and D.~Breitschwerdt, ``{Astrophysical bow shocks:an analytical
  solution for the hypersonic blunt body problem in the intergalactic
  medium.},'' {\em Astronomy and Astrophysics}, vol.~531, A13, 2011.

\bibitem{sfkwbw}
K.~Scherer, F.~Fichtner, J.~Kleimann, T.~Wiengarten, D.~Bomans, and K.~Weis,
  ``{Shock structures of astrospheres.},'' {\em Astronomy and Astrophysics},
  vol.~586, A111, 2016.

\bibitem{schwarze1}
H.~Schwarze, ``{A Uniformly Valid Solution for the Three-Dimensional Inviscid
  Supersonic Flow past Blunt Bodies with Strong Compression in the Shock
  Wave.},'' {\em Acta Mechanica}, vol.~42, pp.~11--35, 1982.

\bibitem{srspkmr4}
S.~R. {Siva Prasad Kochi} and M.~Ramakrishna, ``{Shock capturing with
  discontinuous Galerkin Method using Overset grids for two-dimensional Euler
  equations},'' {\em accepted for publication in Journal of Computational
  Physics}.

\bibitem{paraview}
J.~Ahrens, B.~Geverci, and C.~Law, {\em {ParaView: An End-User Tool for Large
  Data Visualization}}.
\newblock {Elsevier}, 2005.

\bibitem{stYu}
S.-T. Yu, ``{Convenient Method to Convert Two-Dimensional CFD Codes into
  Axisymmetric Ones.},'' {\em Journal of Propulsion and Power}, vol.~9,
  pp.~493--495, 1993.

\bibitem{butcher}
J.~C. Butcher, {\em {Numerical Methods for Ordinary Differential Equations}}.
\newblock {John Wiley and Sons}, 2016.

\bibitem{kxrcf}
L.~Krivodonova, J.~Xin, J.-F. Remacle, N.~Chevaugeon, and J.~Flaherty, ``Shock
  detection and limiting with discontinuous {Galerkin} methods for hyperbolic
  conservation laws.,'' {\em Appl. Numer. Math.}, vol.~48, pp.~323--338, 2004.

\bibitem{srspkmr1}
S.~R. {Siva Prasad Kochi} and M.~Ramakrishna, ``{A compact subcell WENO
  limiting strategy using immediate neighbours for Runge-Kutta discontinuous
  Galerkin methods.},'' {\em International Journal of Computer Mathematics},
  vol.~98, no.~3, pp.~608--626, 2021.

\bibitem{srspkmr3}
S.~R. {Siva Prasad Kochi} and M.~Ramakrishna, ``{A Discontinuous Galerkin
  Overset Scheme Using WENO Reconstruction and Subcells for Two-Dimensional
  Problems},'' {\em Journal of Scientific Computing}, vol.~93:35, 2022.

\bibitem{ks3}
K.~Kitamura and E.~Shima, ``Towards {shock-stable} and accurate hypersonic
  heating computations: {A new pressure flux for AUSM-family schemes}.,'' {\em
  Journal of Computational Physics}, vol.~245, pp.~62--83, 2013.

\bibitem{shu1}
C.-W. Shu, ``Essentially non-oscillatory and weighted essentially
  non-oscillatory schemes for hyperbolic conservation laws.,'' {\em Lecture
  Notes in Mathematics, Springer}, vol.~1697, pp.~325--432, 1998.

\end{thebibliography}

\end{document}